%% file: nystroem_arxiv.tex
\pdfoutput=1
\documentclass[3p,12pt]{elspreprint} 

\input{settings}

\input{commands}

\begin{document}

\input{\CommonPath/title}

%
  \input{\CommonPath/introduction}

%
  \input{\CommonPath/method}

%
  \input{\CommonPath/iga}

%
  \input{\CommonPath/results}

%
  \input{\CommonPath/discussion}

\input{\CommonPath/conclusion}

\clearpage
\appendix

  \input{\CommonPath/appendixTransformation}

%
  \input{\CommonPath/literature}

\end{document}

%% file: settings.tex
\usepackage[utf8]{inputenc} 
\usepackage[english]{babel}  
\usepackage[T1]{fontenc}

\usepackage[intlimits]{amsmath}
\usepackage{amssymb} 
\usepackage{amsthm}
\usepackage{amsfonts}
\usepackage{bm}
\usepackage{mathtools} 

\usepackage{siunitx} 
\usepackage{times}
\usepackage{mathptmx}
\DeclareSymbolFont{cmletters}{OML}{cmr}{m}{n}
\DeclareMathSymbol{\epsilon}{0}{cmletters}{'017} 
\DeclareMathAlphabet{\mathcal}{OMS}{cmsy}{m}{n} 

\usepackage{setspace}
\usepackage{xspace}

\usepackage[noend]{algorithmic}
\usepackage{algorithm}

\usepackage[numbers]{natbib}

\usepackage{comment}
\usepackage{booktabs} 
\usepackage{hyperref} 


%% file: commands.tex
\newcommand{\R}{\mathbb{R}}

\newcommand\op[1]{\mathcal{#1}}
\newcommand\operate[1]{(#1)}

\DeclareMathOperator{\dx}{d}
\DeclareMathOperator{\Tr}{Tr}
\DeclareMathOperator{\dist}{dist}
\DeclareMathOperator{\diam}{diam}
\DeclareMathOperator{\supp}{supp}
\DeclareMathOperator{\diag}{diag}

\renewcommand{\div}{\nabla \cdot}
\newcommand{\grad}{\nabla}
\newcommand{\curl}{\nabla \times}
\newcommand{\dTr}{\op{T}} 
\newcommand{\Laplace}{\Delta}

\newcommand\vek[1]{\mathbf{#1}}
\newcommand\veks[1]{\boldsymbol{#1}}
\newcommand\mat[1]{\mathbf{#1}}
\newcommand\tens[2]{\mathsf{#1}}

\newcommand{\trans}{\intercal} 

\newcommand\pt[1]{\veks{#1}}
\newcommand\ofpt[1]{(\pt{#1})}

\newcommand{\h}{\mathcal{H}} 

\newcommand{\hmatrices}{$\h$-matrices\xspace}

\newcommand\conductivity {k} 
\newcommand\primary      {u} 
\newcommand\dual         {q} 
\newcommand{\BEZ}{B{\'e}zier\xspace}
\newcommand{\NYS}{Nystr{\"o}m\xspace}
\newcommand{\GP}{\pt{y}} 
\newcommand{\locRefInfo}{\bar{\pt{r}}} 
\newcommand{\indexLevel}{\ell}	
\newcommand{\totalLevel}{L}	
\newcommand{\hoelder}{\gamma} 
\newcommand{\Lame}{Lam\'{e}} 

\newcommand{\KV}{\varXi} 	
\newcommand{\KVC}{\varLambda} 	
\newcommand{\Bspline}{B} 	
\newcommand{\CP}{\pt{c}} 		
\newcommand{\CPH}{\CP^h} 		
\newcommand\pu{p}			
\newcommand\puc{\pu_q}			
\newcommand\uu{r} 			
\newcommand\uuc{a} 			
\newcommand\multi{m} 			
\newcommand\cdim{d}			
\newcommand\indexSpan{s} 	
\newcommand\indexA{i} 		
\newcommand\totalA{I} 		
\newcommand\indexB{j} 		
\newcommand\totalB{J} 		
\newcommand\indexC{m} 		
\newcommand\indexD{n} 		
\newcommand{\w}{w} 		
\newcommand{\W}{\w(\uu)} 		
\newcommand{\LshapeCoord}{\zeta} 		

\newcommand\fund[1]{\tens{#1}{2}}

\newcommand\dgamma[1]{\dx s_{\pt{#1}}}

\newcommand{\order}{\mathcal{O}} 

\newcommand\dirac{\delta}
\newcommand\kronecker{\delta}

\newcommand\be{\tau} 
\newcommand\ieg{\tilde{\tau}} 
\newcommand\iegLocal{\grave{\tau}} 
\newcommand{\transMatrix}{\mathbf{T}}
\newcommand{\nodeMatrix}{\mathbf{A}}

\newcommand\Lieg{l} 
\newcommand\LiegLocal{l}

\newcommand\komma{,\,}
\newcommand\und{\,\text{and}\,}

\renewcommand\mit{\,\text{with}\,}
\newcommand\with{\,\text{with}\,}

\newcommand\fromto[2]{\{ #1, \dots, #2 \}}

\newcommand\err{\epsilon}
\newcommand\mychi{\mathcal{X}}

\newcommand{\tikzfig}[5]{
  \begin{figure}[ht]
    \centering
    \includegraphics{#5}
    \caption{#3}
    \label{#4}
  \end{figure}
}

\newcommand{\store}{=} 	
\newcommand{\myinsert}{\leftarrow} 	
 
\newenvironment{myalgorithm}[2]%
{
  \begin{algorithm}
    \caption{#1}
    \label{#2}
    \begin{algorithmic}[1]
}   
{
\end{algorithmic}
\end{algorithm}
}

\newcommand{\todoref}[2][reference]{}
\newcommand{\mytodo}[2][JZBM]{}
\newcommand{\mycomment}[2][-]{}

\newcommand{\mydelete}[2][-]{}
\newcommand\revdel[2]{}
\newcommand\revadd[2]{}
\newcommand\revmod[2]{}
\newcommand\revcomment[1]{} 

\newcommand{\figref}[1]{Figure~\ref{#1}}
\newcommand{\Figref}[1]{Figure~\ref{#1}}
\newcommand{\secref}[1]{section~\ref{#1}}
\newcommand{\Secref}[1]{Section~\ref{#1}}

\providecommand{\algref}[1]{Algorithm~\ref{#1}}
\renewcommand{\algref}[1]{Algorithm~\ref{#1}}


%% file: title.tex
\title{The Isogeometric Nyström Method}

\begin{frontmatter}

\author[ifbaddr]{Jürgen Zechner\corref{cor1}}
\author[ifbaddr]{Benjamin Marussig}
\author[ifbaddr,newcastleaddr]{Gernot Beer}
\author[ifbaddr]{Thomas-Peter Fries}

\address[ifbaddr]{Institute of Structural Analysis, Graz University
  of Technology, Lessingstraße 25/II, 8010 Graz, Austria}

\address[newcastleaddr]{Centre for Geotechnical and Materials Modelling, University of Newcastle,
  Callaghan, NSW 2308, Australia}

\cortext[cor1]{Corresponding author. 
  Tel.: +43 316 873 6181, fax: +43 316 873 6185, mail: \url{ifb@tugraz.at}, web: \url{www.ifb.tugraz.at}}

\input{\CommonPath/abstract}

\end{frontmatter}


%% file: abstract.tex
\begin{abstract}
  In this paper the isogeometric Nyström method is presented. It's outstanding features are: it allows the analysis of domains described by many different geometrical mapping methods in computer aided geometric design and it requires only pointwise function evaluations just like isogeometric collocation methods.
The analysis of the computational domain is carried out by means of boundary integral equations, therefor only the boundary representation is required. The method is thoroughly integrated into the isogeometric framework. For example, the regularization of the arising singular integrals performed with local correction as well as the interpolation of the pointwise existing results are carried out by means of \BEZ elements. 

The presented isogeometric Nyström method is applied to practical problems solved by the Laplace and the \Lame-Navier equation. Numerical tests show higher order convergence in two and three dimensions. It is concluded that the presented approach provides a simple and flexible alternative to currently used methods for solving boundary integral equations, but has some limitations.
\end{abstract}

\begin{keyword}
  Nyström Method \sep 
  Isogeometric Analysis \sep
  Collocation \sep
  Local Refinement \sep 
  Boundary Integral Equation \sep 
\end{keyword}


%% file: introduction.tex
\section{Introduction}
\label{sec:introduction}

Isogeometric analysis has gained an enormous attention in the last decade. It offers a seamless link between the geometrical model to the numerical simulation with no need for meshing \cite{hughes2005}. The key idea is to use the functions describing the geometry in computer aided geometric design (CAGD) also for the analysis.  Since CAGD models are based on boundary representations, a natural integration with simulation is possible using boundary integral equations (BIE).  Most implementations of BIEs discretize the integral equations using basis functions and numerical quadrature to evaluate the bi-linear forms. In this context, the isogeometric boundary element method (BEM) gained much attention recently. In three dimensions, the method has been applied to the Laplace \cite{harbrecht2013}, Stokes \cite{heltai2014}, \Lame-Navier \cite{scott2013,marussig2015}, Maxwell \cite{vazquez2012} and Helmholtz \cite{simpson2014} equations with different discretization methods, i.e. Galerkin \cite{feischl2015} and collocation methods \cite{marussig2015}.

The Nyström method, originally introduced in \cite{nystroem1930}, is an alternative for the numerical solution of BIEs. A unique feature of the method is that the boundary integrals are evaluated directly by means of numerical quadrature without formulating an approximation of the unknown fields. In fact, it is based on \emph{pointwise} evaluations of the fundamental solutions on the boundary of the computational domain. The pointwise nature makes it similar to isogeometric collocation \cite{auricchio2010} that appears to be a very efficient method to solve partial differential equations \cite{schillinger2013}. The Nyström method has been applied to potential \cite{atkinson1997} and electro-magnetic problems \cite{gedney2000,fleming2004}, the Helmholtz equation \cite{canino1998,bruno2001,bremer2014}, Stokes flow \cite{gonzalez2009}, to the analysis of edge cracks \cite{englund2007} and elastic wave scattering \cite{tong2007} and generally, to parabolic BIEs~\cite{tausch2009}.

The implementation requires the treatment of singular fundamental solutions for which several techniques have been developed. Many of them are based on singularity subtraction \cite{anselone1981} or product integration \cite{atkinson1997}. Error analysis for these methods is available \cite{sloan1981} and the interested reader is referred to \cite{hao2014} or textbooks on integral equations such as \cite{atkinson1997} or \cite{kress1999}. In the presented implementation, the locally corrected Nyström method \cite{canino1998} is used for regularization. It is based on the local construction of special quadrature rules for singular functions \cite{strain1995} using composite quadrature rules. An overview of the method can be found in \cite{peterson2009} and \cite{gedney2014}. Mathematical proofs for the solvability and convergence of the locally corrected Nyström method have been provided in \cite{li2013} and \cite{gonzalez2014}. 

Higher order convergence rates can be achieved with the Nyström method based on the order of the chosen quadrature.  A requirement for higher order convergence is the continuity of the boundary representation which has to be smooth \cite{atkinson1997}. As a consequence, standard triangulations are insufficient for curved geometries due to the potential kink between patches along their common edge. With the help of re-parametrization, patches with the required continuity may be constructed \cite{ying2006} but the procedure is complex. This is why the Nyström method has so far been mostly applied to analytical surfaces. Its application to CAGD geometry descriptions has the advantage that continuity of the boundary can be better controlled, as already noted in \cite{canino1998}. However, real-world geometries still contain corners and edges by design. In order to restore higher order convergence, the elements which define the integration regions are graded \cite{atkinson1997}. For two-dimensional problems, several other approaches exist \cite{kress1990,helsing2008,bruno2009,bremer2010,gillman2014}. While it is claimed there that these methods can be extended to three-dimensional problems as well, only few results are reported \cite{bremer2012a}. A solution for tackling singularities, arising at regions with mixed boundary conditions is presented in \cite{cheng1994} and in \cite{helsing2009}.

In the following sections, the \emph{isogeometric Nyström method} is introduced for the Laplace and the \Lame-Navier equation. The formulation is based on CAGD boundary representations which are partitioned into integration regions for the composite quadrature. Therefor, the arbitrary selectable CAGD technology is only required to provide a valid geometrical mapping from the parameter to the real space. The numerical scheme consists of point collocation on the surface of the computational domain.  
The regularization by means of local correction is performed with \BEZ elements. To preserve higher order convergence, a priori grading of elements at corners and edges is realized in the parameter space. Without loss of generality, the procedure is adapted to boundary representations based on non-uniform rational B-splines (NURBS). For tensor product descriptions of surfaces the authors present a strategy for the \emph{local refinement} of elements. For post-processing purposes, the pointwise existing results are interpolated over the boundary again with the help of \BEZ elements.

In summary, the analysis with the proposed isogemetric Nyström method offers the following main advantages to both, the IGA and the Nyström community:
\begin{itemize}
\item isogeometric collocation scheme for boundary integral equations,
\item free choice of CAGD technology taken for the boundary representation,
\item local refinement even for tensor product patches,
\item application of the Nyström method to real-world geometries.
\end{itemize}

The paper is organized in the following parts: In \secref{sec:bie} the boundary integral equation for heat conduction and elasticity is revisited.  The principles and requirements of the locally corrected Nyström method are explained.  \Secref{sec:iga} describes the application and adaptation of the isogeometric framework to the method.  Emphasis is given to NURBS surface descriptions and to the local refinement.  In \secref{sec:results} several numerical tests in two and three dimensions are presented.  These results as well as the advantages and the limitations of the isogeometric Nyström method are discussed in \secref{sec:discussion}.  The paper closes with concluding remarks in \secref{sec:conclusion}.


%% file: method.tex
\section{Boundary Integral Equations}
\label{sec:bie}

\mytodo[JZ]{%
  \begin{itemize}
  \item [\done] Local Correction
  \item Erklärung, was ist ``Order of quadrature'' und die Konvergenz!
  \end{itemize}
}

Mathematical models described by mixed elliptic boundary value problems (BVP) are considered in this section as well as their discretization with the Nyström method. They are based on Laplace's equation and on the \Lame-Navier equation which describe isotropic steady heat conduction and linear elasticity respectively. 

\subsection{Boundary Value Problem}
\label{sec:bie:bvp}

Let $\Omega$ be the considered domain, $\primary$ a generalized unknown and $\op{L}$ an elliptic partial differential operator. 
Then the BVP is generalized to the following form: Find $\primary\ofpt{x}$ so that
\begin{equation}
\label{eq:BVP}
  \begin{aligned}
    \op{L} \primary \ofpt{x} & = 0 
    & & \forall \pt{x} \in \Omega\\
    \dTr \primary \ofpt{x} & = \dual \ofpt{y}  = {g}_N\ofpt{y}
    & & \forall \pt{y} \in \Gamma_N\\
    \Tr \primary \ofpt{x} &= \dual \ofpt{y} = g_D\ofpt{y} 
    & & \forall \pt{y} \in \Gamma_D.
  \end{aligned}
\end{equation}
For \eqref{eq:BVP} the boundary $\Gamma=\partial\Omega$ is split into a Neumann and a Dirichlet part ${\Gamma = \Gamma_N \cup \Gamma_D}$ so that ${\Gamma_N \cap \Gamma_D = \emptyset}$. 
The boundary trace 
\begin{align}
  \label{eq:trace}
  \Tr \primary \ofpt{x} &= \lim_{\vek{x}\rightarrow\vek{y}} \primary \ofpt{x} =
  \primary \ofpt{y} & & \pt{x} \in \Omega,\, \pt{y} \in \Gamma
\end{align}
maps the primal field~$\primary\ofpt{x}$ in the domain $\Omega$ to $\primary\ofpt{y}$ on the boundary $\Gamma$. 
The prescribed Neumann and Dirichlet data are denoted by ${g}_N$ and ${g}_D$ respectively.

For the heat equation the partial differential operator is defined by the Laplacian $\Laplace=\div \grad$
\begin{equation}
  \label{eq:LaplaceEquation}
  \op{L} \primary\ofpt{x} := -\conductivity \Laplace \primary\ofpt{x}
\end{equation}
where $k$ denotes the conductivity and $\primary$ the temperature. 
The conormal-derivative maps temperature $\primary$ to heat flux $\dual$ on the surface $\Gamma$ and is defined by the normal derivative
\begin{align}
  \label{eq:LaplaceEquationConormal}
  \dTr \primary \ofpt{x} & = %
  \conductivity \grad\primary\ofpt{y} \cdot \vek{n}\ofpt{y} = \dual\ofpt{y}  & & %
  \pt{x} \in \Omega \komma \pt{y} \in
  \Gamma
\end{align}
with $\vek{n}$ denoting the unit outward normal vector of $\Omega$.

For the \Lame-Navier equation
\begin{equation}
  \label{eq:LameEquation}
  \op{L} \primary\ofpt{x} := -\left( \lambda + 2 \mu \right) 
  \div \grad \primary \ofpt{x} + \mu \curl (\curl \primary \ofpt{x}) =
  0
\end{equation}
holds with $\primary$ denoting the displacement and with the \Lame-constants $\lambda$ and $\mu$ \cite{kupradze1979}. For elasticity the conormal derivative is
\begin{align}
  \label{eq:LameEquationConormal}
  \dTr \primary \ofpt{x} & = %
  \lambda \div \primary \ofpt{y} \vek{n}\ofpt{y} + %
  2\mu \grad \primary \ofpt{y} \cdot \vek{n}\ofpt{y} + %
  \mu \vek{n}\ofpt{y} \times (\curl \primary \ofpt{y}) = %
  \dual\ofpt{y}& & %
  \pt{x} \in \Omega,\, \pt{y} \in \Gamma
\end{align}
which maps displacements $\primary$ to surface traction $\dual$.

\subsection{Boundary Integral Equations}
\label{sec:bie:bie}

The variational form of the BVP \eqref{eq:BVP} can be solved by means of boundary integral equations \cite{atkinson1997,beer2008,sauter2011}. For the solution of many physical problems, Fredholm integral equations of the first kind
\begin{align}
  \label{eq:FredholmFirst}
  \primary\ofpt{x} & = \operate{\op{V} \phi} \ofpt{x} & & \forall \pt{x} \in \Gamma
\end{align}
or the second kind
\begin{align}
  \label{eq:FredholmSecond}
  \primary\ofpt{x} & = \operate{(\op{C} + \op{K}) \psi }\ofpt{x} & & \forall \pt{x} \in \Gamma
\end{align}
are used. In these equations
\begin{align}
  \label{eq:SL}
  \operate{\op{V} \dual }\ofpt{x} &= \int_\Gamma \fund{U}(\pt{x},\pt{y}) \dual \ofpt{y}
  \dgamma{y} & & \forall \pt{x},\pt{y} \in \Gamma
\end{align}
denotes the single layer and
\begin{align}
  \label{eq:DL}
  \operate{\op{K} \primary}\ofpt{x} &= %
  \int_{\Gamma} \dTr\fund{U}(\pt{x},\pt{y}) \primary \ofpt{y}
  \dgamma{y} & & %
  \forall \pt{x},\pt{y} \in \Gamma 
\end{align}
the double layer boundary integral operator. 
The \emph{kernel function} $\fund{U}$ is the fundamental solution for the underlying problem which depends on the Euclidean distance {$r=|\pt{x}-\pt{y}|$}. 
In case of Laplace's equation it is
\begin{equation}
  \label{eq:LaplaceFundSol}
  \fund{U}(\pt{x},\pt{y}) = \frac{1}{4\pi r}
\end{equation}
for $\Omega \in \R^3$. In case of the \Lame-Navier equation it is Kelvin's fundamental solution \cite{beer2008}
\begin{equation}
  \label{eq:LameFundSol}
  \fund{U}_{ij} (\pt{x},\pt{y}) = %
  \frac{1}{8\pi r} \frac{\lambda+\mu}{\mu(\lambda+2\mu)} %
  \left( %
    r_{,i}r_{,k} + \frac{\lambda + 3\mu}{\lambda + \mu} \dirac_{ij} %
  \right) 
\end{equation}
represented as a tensor with $i, j=\fromto{1}{3}$. The fundamental solutions tend to infinity if {$r\rightarrow 0$}. Therefore the integrals in \eqref{eq:SL} and \eqref{eq:DL} are and need to be regularized or evaluated analytically \cite{sauter2011}. For \eqref{eq:DL} the regularization results in an integral free jump term 
\begin{align}
  \label{eq:JumpTerm}
  \operate{\op{C}  \primary} \ofpt{x} && \with && \op{C}=\tfrac{1}{2}\op{I}
\end{align}
on smooth surfaces.

For indirect formulations with first \eqref{eq:FredholmFirst} or second kind integral equations \eqref{eq:FredholmSecond}, the unknowns $\phi$ and $\psi$ are usually non-physical quantities and only intermediate results in order to evaluate quantities in the interior $\pt{x} \in \Omega$ by means of the representation formulas

\begin{align}
  \label{eq:RepresentationFormulas}
  \primary\ofpt{x} & = %
  \operate{\op{V} \phi} \ofpt{x} & & \und & &%
  \primary\ofpt{x} = %
  \operate{\op{K} \psi} \ofpt{x} & & %
  \forall \pt{x} \in \Omega.
\end{align}
Working with physical quantities $\primary$ and $\dual$ on the boundary $\Gamma$ only, the BVP \eqref{eq:BVP} is solved by means of a direct boundary integral formulation
\begin{align}
  \label{eq:DirectBIE}
  \operate{(\op{C} + \op{K}) \primary}\ofpt{x} & = %
  \operate{\op{V} \dual} \ofpt{x} & & %
  \forall \pt{x} \in \Gamma
\end{align}
instead.

\subsection{Discretization with the Nyström Method}
\label{sec:bie:nystroem}

In his original paper, \citet{nystroem1930} proposed the discretization of second kind integral equations \eqref{eq:FredholmSecond} by means of a numerical quadrature
\begin{align}
  \label{eq:NystroemDiscretization}
  \operate{(\tfrac{1}{2}\op{I} + \op{K}) \psi } (\pt{x}_i) & \approx %
  c\psi(\pt{x}_i) + \sum_{j=1}^{n}
  \dTr\fund{U}(\pt{x}_i,\pt{y}_j)\psi(\pt{y}_j) w_j && %
  \pt{x}_i \in \Gamma \komma \pt{y}_j \in \Gamma .
\end{align}
In this equation, $\GP_j$ are the evaluation points of the $n$-point numerical quadrature and $w_j$ are their corresponding weights. In order to set up a linear system of equations, the quadrature sum \eqref{eq:NystroemDiscretization} is collocated at distinct points $\pt{x}_i$ with $i=\fromto{1}{n}$ resulting in the matrix equation
\begin{equation}
  \label{eq:IndirectMatrixDLP}
  \left( \tfrac{1}{2}\mat{I} + \mat{K} \right) \vek{\psi} = \vek{u}
\end{equation}
on smooth surfaces. For Fredholm integral equations of the first kind, the system is
\begin{equation}
  \label{eq:IndirectMatrixSLP}
  \mat{V} \vek{\phi} = \vek{u}.
\end{equation}
For a direct formulation with mixed boundary conditions, a block system of equations
\begin{align}
  \label{eq:DirectMatrixEquation}
  \begin{matrix}
    \pt{x}\in \Gamma_{D}: \\
    \pt{x}\in \Gamma_{N}: 
  \end{matrix} \quad
  \begin{pmatrix}
    \vek{V}_{DD} & -{\vek{K}}_{DN} \\
    \vek{V}_{ND} & -{\vek{K}}_{NN}
  \end{pmatrix}
  \begin{pmatrix}
    \vek{\dual}_{D} \\ 
    \vek{\primary}_{N}
  \end{pmatrix}
  = 
  \begin{pmatrix}
    {\vek{K}}_{DD} & -\vek{V}_{DN} \\
    {\vek{K}}_{ND} & -\vek{V}_{NN}
  \end{pmatrix}
  \begin{pmatrix}
    \vek{g}_D \\ 
    \vek{g}_N
  \end{pmatrix}
  \begin{matrix} 
    \\
    .
  \end{matrix} 
\end{align}
is taken like the formulation presented in \cite{zechner2013}. If integral free jump terms are present, they are already integrated in the system matrices $\mat{K}$. If the surface $\Gamma$ is smooth and the kernel function $\fund{U}$ is regular, entries of the system matrix only consist of pointwise evaluations 
\begin{align}
  \label{eq:PointEvaluation}
  \mat{V}[i,j] & = \fund{U}(\pt{x}_i,\pt{y}_j) w_j & & \und & & %
  \mat{K}[i,j] = \kronecker_{ij}c+\dTr\fund{U}(\pt{x}_i,\pt{y}_j) w_j.
\end{align}
For the considered applications the kernel functions are singular with $c=\tfrac{1}{2}$ on smooth surfaces. Moreover, the fundamental solution is undefined if $\pt{x}_i=\pt{y}_j$ so that special treatment is necessary to evaluate the corresponding system matrix entries. 

Based on a technique for the construction of quadrature rules with arbitrary order for given singular functions presented in \cite{strain1995}, the authors of \cite{canino1998} developed the locally corrected Nyström method for the solution of the Helmholtz equation. This particular regularization technique is taken for the framework presented in this paper. The main idea is to replace the contribution of the original kernel function in the neighborhood $\Omega_x$ of the collocation point $\pt{x}_j$ with a corrected regular one, so that the new kernel function is defined by
\begin{equation}
  \label{eq:PointEvaluationCorrected}
  \fund{U^{*}}(\pt{x}_i,\pt{y}_j) = %
  \begin{cases} 
    \fund{L}(\pt{x}_i,\pt{y}_j) & \textnormal{if } \pt{x}_i \in \Omega_x \\
    \fund{U}(\pt{x}_i,\pt{y}_j) & \textnormal{otherwise.} \end{cases}
\end{equation}
The locally corrected kernel $\fund{L}$ for the collocation point $\pt{x}_i$ is computed at $n$ corresponding field points $\pt{y}_j \in \Omega_x$ by solving the linear system
\begin{align}
  \label{eq:LocalCorrection}
  \sum_{j=1}^n N_i(\pt{y}_j) L(\pt{x}_i,\pt{y}_j) w_j & =%
  \int_{\Gamma \cap \Omega_x} \fund{U}(\pt{x}_i,\pt{y}) N_i(\pt{y}) \dgamma{y} && \with && i=\fromto{1}{m}.
\end{align}
Equation \eqref{eq:LocalCorrection} introduces a space of $m$ testfunctions $N_i$, hence the singularity of the original kernel is treated in a weak sense on the right hand side. The choice of the test functions in the presented application is discussed in \secref{sec:iga}. Because of the singularity of $\fund{U}$, treating the moments on the right hand side of \eqref{eq:LocalCorrection} requires regularization. This equation constructs a numerical quadrature, where the weights $w_j$ are not explicitly calculated but collected together with the corrected kernel to $\tilde{w}_j=L(\pt{x}_i,\pt{y}_j) w_j$. Finally, the linear system in matrix form is
\begin{align}
  \label{eq:LocalCorrectionSystem}
  \mat{N} \tilde{\vek{w}} & = \vek{g} & & %
  \mit & & %
  \mat{N} \in \R^{m\times n}\komma %
  \tilde{\vek{w}} \in \R^n \komma %
  \vek{g} \in \R^m .
\end{align}
The matrix $\mat{N}$ consists of evaluations of the test-functions $N_i$ and $\vek{g}$ contains accurately evaluated singular moments.  Equation \eqref{eq:LocalCorrectionSystem} is numerically solved for $\tilde{\vek{w}}$ with $LU$-decomposition if $n=m$.  The number $m$ of chosen test functions $N_i$ may be smaller or larger than the number $n$ of sample points $\pt{y}_i$.  In that case, a valid solution is found by means of least-squares or a minimum norm solution \cite{strain1995}.  

In this paper, the order of a quadrature $p$ is defined as the degree of the highest polynomial that it does integrate exactly. Therefor, the numerical solution of equations \eqref{eq:IndirectMatrixDLP}, \eqref{eq:IndirectMatrixSLP} or \eqref{eq:DirectMatrixEquation} converge with $p+1=n+1$ for an $n$-point quadrature. Although composite Gauss quadrature is used in the presented formulation, the order of convergence for the method does not reach $2n+1$ due to the local correction \cite{strain1995}.

In practical applications, the considered surface $\Gamma$ of the domain is approximated by $\Gamma_h$ consisting of non-overlapping patches $\be$ so that
\begin{align}
  \Gamma_h = \bigcup_{l=1}^{L} \be_l.
\end{align}
As a requirement for the locally corrected Nyström method, a composite quadrature based on $\Gamma_h$ is chosen for the discretization of the integral equation. Additionally, the quadrature is chosen to be of open type, which means that no quadrature points are located at the boundary of the integration region. This is because such points can be located at physical edges, where the integral kernel may be undefined or diverging. This particular choice comes with an advantage: In contrast to the BEM and due to its pointwise nature, the Nyström method does not require or impose any connectivities between patches. Hence, non-conforming meshes are supported inherently.

To converge with respect to the order of the underlying quadrature, the Nyström method for solving integral equations of the kind \eqref{eq:FredholmFirst}, \eqref{eq:FredholmSecond} or \eqref{eq:DirectBIE} requires a smooth surface. Usually, this is not feasible by means of a standard triangulation of $\Gamma$.  Therefor, the authors propose the application of CAGD surface descriptions, which fulfill this requirement.  Moreover, $\Gamma_h=\Gamma$ which means that the discretization introduces no geometry error.


%% file: iga.tex
\section{Isogeometric Framework}
\label{sec:iga}

In this section the isogeometric paradigm is combined with the locally corrected \NYS method.  The term Cauchy data is taken to refer to quantities on the boundary $\Gamma$ appearing in the discrete boundary integral equations~\eqref{eq:IndirectMatrixDLP}, \eqref{eq:IndirectMatrixSLP} and~\eqref{eq:DirectMatrixEquation}.

The key aspect of the isogeometric concept is to utilize the boundary representation of design models directly in the analysis.  Thus, it has been applied to a variety of models with different surface descriptions such as subdivision surfaces \cite{cirak2000}, tensor product surfaces \cite{hughes2005} and T-spline surfaces \cite{scott2013}. Since the most commonly used CAGD technology in engineering design are tensor product surfaces and based on non-rational B-splines (NURBS), the paper focuses on geometry descriptions based on this approach. However, the implementation of the \NYS method to other surface representations is straight forward, since the approximation of the Cauchy data as well as the partitioning of the elements for the integration is independent of the geometrical parametrization.  In fact, the only requirement is a valid geometrical mapping $\mychi\ofpt{\uu}$ from local coordinates $\vek{\uu}=(\uu_1,\dots,\uu_{\cdim-1})^\trans$ to global coordinates $\pt{x}=(x_1,\dots,x_{\cdim})^\trans$ in the $d$-dimensional Cartesian system $\R^\cdim$.  To be precise the Gram's determinant has to be non-singular. The corresponding Gram matrix is given by
\begin{align}
  \label{eq:GramMatrix}
  \mat{G}\ofpt{\uu} &\coloneqq \mat{J}_\mychi^\trans \ofpt{\uu} \: \mat{J}_\mychi \ofpt{\uu} 
  \in \R^{ \left( \cdim-1 \right) \times \left( \cdim-1 \right)}
\end{align}
where the Jacobi-matrix 
\begin{align}
  \label{eq:JacobiMatrix}
   \mat{J}_\mychi \ofpt{\uu} &\coloneqq \left(  \frac{\partial \mychi_\indexA}{ \partial \uu_\indexB}  \right)%
   && \mit && \indexA = \fromto{1}{\cdim}\komma \indexB = \fromto{1}{\cdim-1}
   \qquad 
\end{align}
results form the geometrical mapping $\mychi\ofpt{\uu}$~\cite{sauter2011}.

\subsection{Geometry Representation}
\label{sec:iga:geometry}

CAGD objects are usually defined by a set of boundary curves or surfaces. 
They are generally referred to as \emph{patches} for the remainder of this paper.
A significant property is that the continuity within patches can be controlled by the associated basis functions.
Hence, such patches are able to represent a smooth geometry without any artificial corners or edges. 

\subsubsection{Basis functions}

The \emph{knot vector}~$\KV$ is the fundamental element for the construction of the basis functions.  It is characterized as a non-decreasing sequence of coordinates~{$\uu_\indexA \leq \uu_{\indexA+1}$} which defines the parametric space of a patch.
The coordinates itself are called \emph{knots} and the half-open interval $\left[\uu_i, \uu_{i+1}\right)$ is called \emph{knot span}. Knot values may not be unique which is referred to as the \emph{multiplicity} of a knot, which is then larger than one. Together with a corresponding polynomial degree $\pu$ a set of piecewise polynomial basis functions $\Bspline_{\indexA,\pu}$, so called \emph{B-splines} are defined recursively. After introducing piecewise constant functions ($\pu=0$)
\begin{equation}
  \label{eq:Bspline_N0}
  \Bspline_{\indexA,0}(\uu) = 
  \begin{cases} 
    1 & \textnormal{if } \uu_{\indexA}\leqslant \uu < \uu_{\indexA+1}\\
    0 & \textnormal{otherwise,}
  \end{cases}
\end{equation}
higher degree B-splines are constructed as a strictly convex combination of basis functions of the previous degree
\begin{equation}
  \label{eq:Bspline_Np}
  \Bspline_{\indexA,\pu}(\uu)  = \frac{\uu-\uu_{\indexA}}{\uu_{\indexA+\pu}-\uu_{\indexA}} \: \Bspline_{\indexA,\pu-1}(\uu) 
  + \frac{\uu_{\indexA+p+1}-\uu}{\uu_{\indexA+p+1}-\uu_{\indexA+1}} \: \Bspline_{\indexA+1,\pu-1}(\uu).
\end{equation}
Further, their first derivatives are also a linear combination of B-splines of the previous degree
\begin{equation}
	\label{eq:Bspline_Np_Deriv}
	\frac{\partial \Bspline_{\indexA,\pu}(\uu) }{\partial \uu} = \Bspline^\prime_{\indexA,\pu}(\uu) 
	= \frac{\pu}{\uu_{\indexA+\pu}-\uu_{\indexA}} \: \Bspline_{\indexA,\pu-1}(\uu) 
	- \frac{\pu}{\uu_{\indexA+\pu+1}-\uu_{\indexA+1}} \: \Bspline_{\indexA+1,\pu-1}(\uu).
\end{equation}
The support $\supp{ \{\Bspline_{\indexA,\pu} \} }=\fromto{\uu_{\indexA}}{\uu_{\indexA+\pu+1}}$ is local and entirely defined by $\pu + 2$ knots.  $\Bspline_{\indexA,\pu}$ is described by a polynomial segment within each non-zero knot span $\left(\uu_{\indexSpan} < \uu_{\indexSpan+1} \right)$ of its support.  The continuity between adjacent segments is $C^{\pu-\multi}$ where $\multi$ denotes the multiplicity of the joint knot. Consequently, the continuity of the basis functions $\Bspline_{\indexA,\pu}$ at their knots is determined by the corresponding knot vector.  \Figref{fig:bspline} illustrates this relation for two quadratic B-splines. Note that the continuity between the polynomial segments decreases at the double knot.
\tikzfig{tikz/bsplineSegments}{1.1}{Examples of two quadratic B-splines ($\pu=2$) with different continuities at their middle knots. The knot vector of the left B-spline is defined by $\KV = \left\{1,2,3,4\right\}$, $\KV = \left\{1,2.5,2.5,4\right\}$ defines the knot vector of the right one. The resulting polynomial segments are indicated by dashed lines, solid lines represents the B-splines.}{fig:bspline}{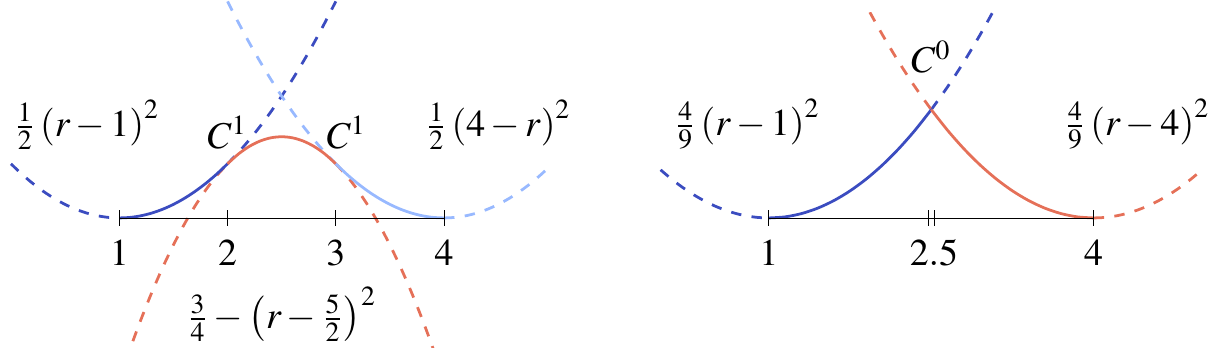}
In general, the knots are arbitrarily distributed. This is emphasized by referring to the knot vector as \emph{non-uniform}.  The term \emph{open} knot vector indicates that the first and last knot are $C^0$-continuous.  A special knot sequence is
\begin{align}
  \label{eq:knot_vector_BEZ}
  \KV = \left\{ \uu_0 = \dots = \uu_\pu, \uu_{\pu+1} = \dots = \uu_{2\pu+1} \right\}
\end{align}
where the multiplicity of all knots is equal to the polynomial order, i.e. $\pu + 1$. 
The resulting basis functions $\Bspline_{\indexA,\pu}$ are classical $\pu$th-degree Bernstein polynomials which extend over a single non-zero knot span. 

\subsubsection{Curves}
The geometrical mapping of B-spline curves of degree $\pu$ is
\begin{align}
  \label{eq:Bspline_mapping}
  \mychi\operate{\uu} &\coloneqq \pt{x}\operate{\uu} = %
  \sum_{\indexA=0}^{\totalA-1} \Bspline_{\indexA,\pu}\operate{\uu} \: \CP_{\indexA} 
\end{align}
and its Jacobi-matrix is given by
\begin{align}
  \label{eq:Bspline_mapping_jac}
  \mat{J}_\mychi \operate{\uu} &\coloneqq %
  \sum_{\indexA=0}^{\totalA-1} \Bspline^\prime_{\indexA,\pu}\operate{\uu} \: \CP_{\indexA}
\end{align}
where $\totalA$ is the total number of basis functions due to a knot vector $\KV_\totalA$ and the \emph{control points} $\CP_{\indexA}$ are the corresponding coefficients in physical space. The resulting piecewise polynomial curve $\pt{x}\operate{\uu}$ inherits the continuity properties of its underlying basis functions. If they span only over a single non-zero knot span, i.e. $\KV_\totalA$ is of form \eqref{eq:knot_vector_BEZ}, the curve is referred to as \BEZ curve.

In order to represent models like conic sections which are based on rational functions, \emph{weights} $\w_\indexA$ are introduced. They are associated with the control points yielding to \emph{homogeneous coordinates} given by 
\begin{align}
  \CPH_\indexA & = (\w_\indexA \CP_\indexA, \w_\indexA)^\trans = (\CP^\w_\indexA, \w_\indexA)^\trans\in \R^{\cdim+1}.
\end{align}
Applying the mapping~\eqref{eq:Bspline_mapping} to $\CPH_\indexA$ defines a B-spline curve $\pt{x}^h\operate{\uu}$ in the projective space $\R^{\cdim+1}$. Hence, the geometrical mapping has to be extended by a perspective mapping $\op{Y}$ with the center at the origin of $\R^{\cdim+1}$.%
\tikzfig{tikz/homogeneousCoordinates}{2}{Perspective mapping $\op{Y}$ of a quadratic B-spline $\pt{x}^h\operate{\uu}$ in homogeneous form $\R^{3}$ to a circular arc $\pt{x}\operate{\uu}$ in physical space $\R^{2}$.}{fig:projection}{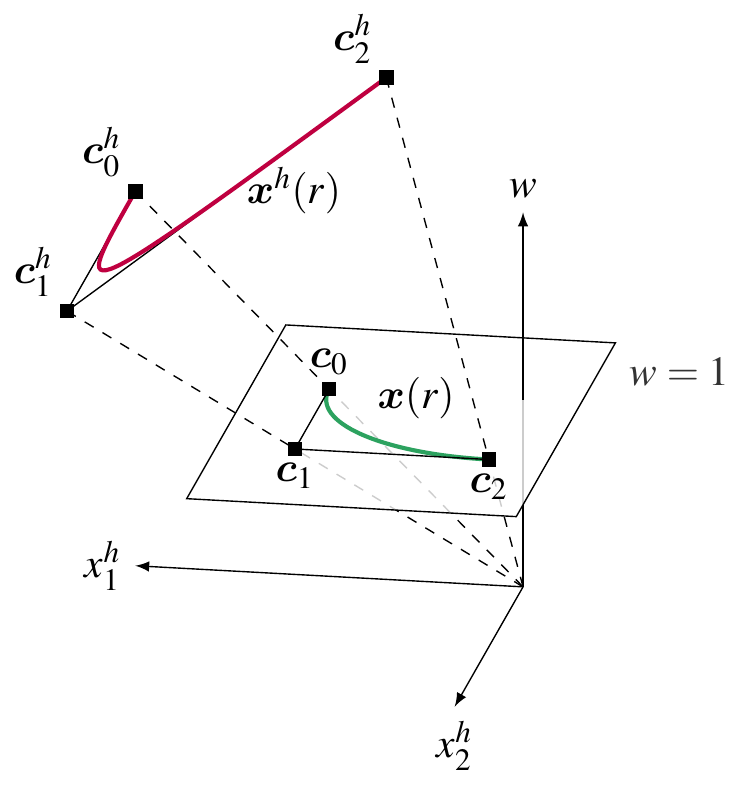}%
This is depicted in \figref{fig:projection} for a circular arc in $\R^2$. Its homogeneous representation is defined by quadratic B-splines and the control points $\CPH_0 =(2,0,2)^\trans$, $\CPH_1 =(2,3,2)^\trans$ and $\CPH_2 =(0,4,4)^\trans$ in the projective space $\R^3$. The homogeneous vector components {$\pt{x}^\w = \left(\pt{x}^h_1,\dots,\pt{x}^h_\cdim \right)^\trans$} are mapped to the physical space by
\begin{align}
  \label{eq:NURBS_curve_mapping}
  \pt{x}\operate{\uu} &= \op{Y}\operate{\pt{x}^h\operate{\uu}} =  \frac{\pt{x}^\w\operate{\uu} }{\w\operate{\uu}}
  && \with && \w\operate{\uu} = \sum_{\indexA=0}^{\totalA-1} \Bspline_{\indexA,\pu} (\uu) \w_{\indexA}.
\end{align}
The resulting rational function which describes such a curve is called NURBS (non-uniform rational B-spline).  The geometrical mapping with NURBS is defined by
\begin{align}
  \label{eq:NURBS_curve_mapping_jac}
  \mat{J}_\mychi \operate{\uu} &\coloneqq %
  \frac{ \W \frac{\partial \pt{x}^\w\operate{\uu} }{ \partial \uu} - \frac{\partial \W }{ \partial \uu} \pt{x}^{\w}(\uu)  }{\left(\W\right)^2}
\end{align}
where the derivatives are defined by
\begin{align}
  \frac{\partial \W }{ \partial \uu}  &= \sum_{\indexA=0}^{\totalA-1} \Bspline^\prime_{\indexA,\pu} (\uu) \w_{\indexA}
  && \und &&
  \frac{\partial \pt{x}^\w\operate{\uu} }{ \partial \uu} = \sum_{\indexA=0}^{\totalA-1} \Bspline^\prime_{\indexA,\pu} (\uu) \CP^\w_\indexA.
\end{align}	

\subsubsection{Tensor Product Surfaces}

B-spline surfaces are defined by tensor products of univariate basis functions which are related to separate knot vectors $\KV_\totalA$ and $\KV_\totalB$. 
The geometrical mapping is given by 
\begin{align}
	\label{eq:Bspline_mapping_patch}
	\mychi\ofpt{\vek{\uu}} &\coloneqq \pt{x}\ofpt{\vek{\uu}} = 
	\sum_{\indexA=0}^{\totalA-1} \sum_{\indexB=0}^{\totalB-1} \Bspline_{\indexA,\pu_1} \operate{\uu_1}  \Bspline_{\indexB,\pu_2}\operate{\uu_2}    \: \CP_{\indexA,\indexB}
\end{align}
with $\pu_1$ and $\pu_2$ denoting different polynomial degrees.
The corresponding Jacobian is calculated by substituting the basis functions by its derivatives, alternately for each parametric direction.
Further, a surface is labeled as \BEZ patch, if $\KV_\totalA$ and $\KV_\totalB$ are of form \eqref{eq:knot_vector_BEZ} and the derivation of NURBS surfaces is analogous to NURBS curves.

In general the surface representation by tensor products is an extremely efficient technique compared to other geometry descriptions \cite{de-boor1978}.

However, the applicability of the tensor product approach is limited. Especially in the context of refinement which can not be performed locally. In the scope of the present work, this drawback does not limit the ability for local refinement with the \NYS method.

\subsection{Approximation of Cauchy Data}
\label{sec:iga:cauchy}

\subsubsection{Element-wise Discretization}
In order to evaluate the boundary integral properly, the patch $\be$ is subdivided into a set of \emph{elements} $\ieg$, so that
\begin{align}
  \be = \bigcup_{\indexA=1}^{\totalA} \ieg_\indexA.
\end{align}
A distinguishing feature of the \NYS method is that the element-wise discretization of the Cauchy data is directly expressed through points defined by the applied quadrature rule on the geometry, rather than control points. The quadrature points are distributed with respect to their coordinates $\xi \in \left[-1,1 \right]$ 
in the reference element. Consequently, the \emph{$p$-refinement strategy} is determined by the increase of the quadrature order and thus, the quadrature points on the element. To resolve their location $\GP$ on $\ieg$ in physical space, the mappings $\mychi_{\xi}(\xi)$ for the reference element and $\mychi(\uu)$ for the geometry are consecutively applied. For a curved element $\be=\ieg$ these mappings are depicted in \figref{fig:GaussPoints}.%
%
%
%
\tikzfig{tikz/quadratureDistributionSingleElement}{1.2}{The distribution of the quadrature points $\GP$ in the reference space $\xi$, the parameter space $\uu$ of the geometry and the physical space $\pt{x}\in\R^2$. The quadrature points are indicated by circles. Squares denote control points $\CP$.}{fig:GaussPoints}{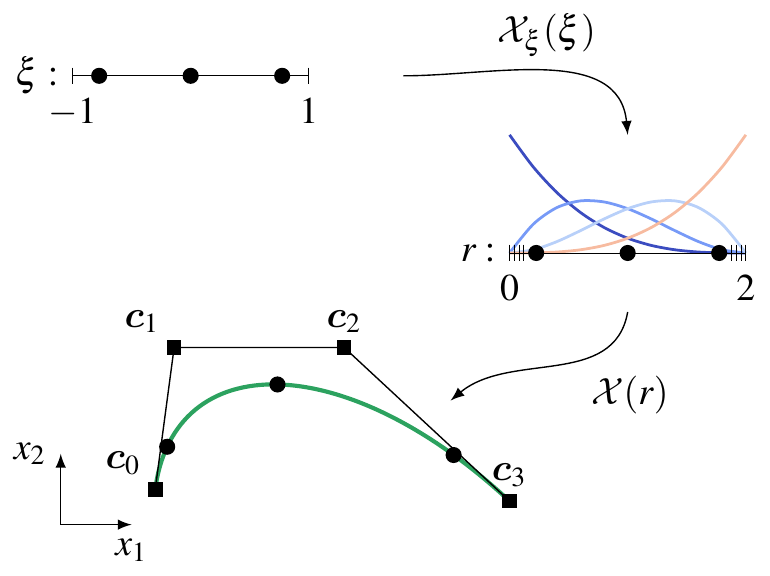}
%

In contrast to other isogeometric methods, the Cauchy data is not expressed in terms of variables at control points $\CP$. For the analysis with the isogeometric Nyström method, Cauchy data are taken from evaluations in the quadrature points $\GP\in\Gamma$ directly, without any transformation.

\subsubsection{Patch-wise Discretization}
Generally, the arrangement of $\ieg$ over $\be$ is independent of the geometry. But its representation has to be smooth within each $\ieg$. As indicated in \secref{sec:iga:geometry}, the continuity of the geometry is directly linked to the knot vector $\KV$. It is convenient, to define $\ieg$ by means of an artificial knot vector $\KVC$ so that $\KV \subset \KVC$. To be precise, the purpose of $\KVC$ is not to construct basis functions but to organize the global partition of $\ieg$ properly. Consequently, $\uuc$ denote knots of $\KVC$. Based on the initial discretization, i.e. $\KV = \KVC$, the approximation quality of the Cauchy data can be improved by inserting additional \emph{unique} knots $\bar{\uuc_\indexA}$ into a knot span $\indexSpan$ of $\KVC$, so that $\uuc_{\indexSpan} < \bar{\uuc_\indexA} < \uuc_{\indexSpan+1}$.  This procedure defines an \emph{$h$-refinement strategy} that is performed in parametric space and preserves the continuity requirements with respect to the geometry. As indicated in \figref{fig:hRef} this is independent of the geometry representation.
\tikzfig{tikz/quadratureDistribution}{1.2}{The geometry of a cubic B-spline curve is defined by the knot vector $\KV = \left\{0,0,0,0,2,4,4,4,4\right\}$ and the control points $\CP_\indexA$ with $\indexA=\fromto{0}{4}$. The elements $\ieg_j$ with $j=\fromto{1}{4}$ along the curve are defined by $\KVC = \left\{0,0,0,0,1,2,3,4,4,4,4\right\}$ in the parameter space. Each element is equipped with a two point quadrature rule defined in the reference interval $\xi \in \left[-1,1 \right]$.}{fig:hRef}{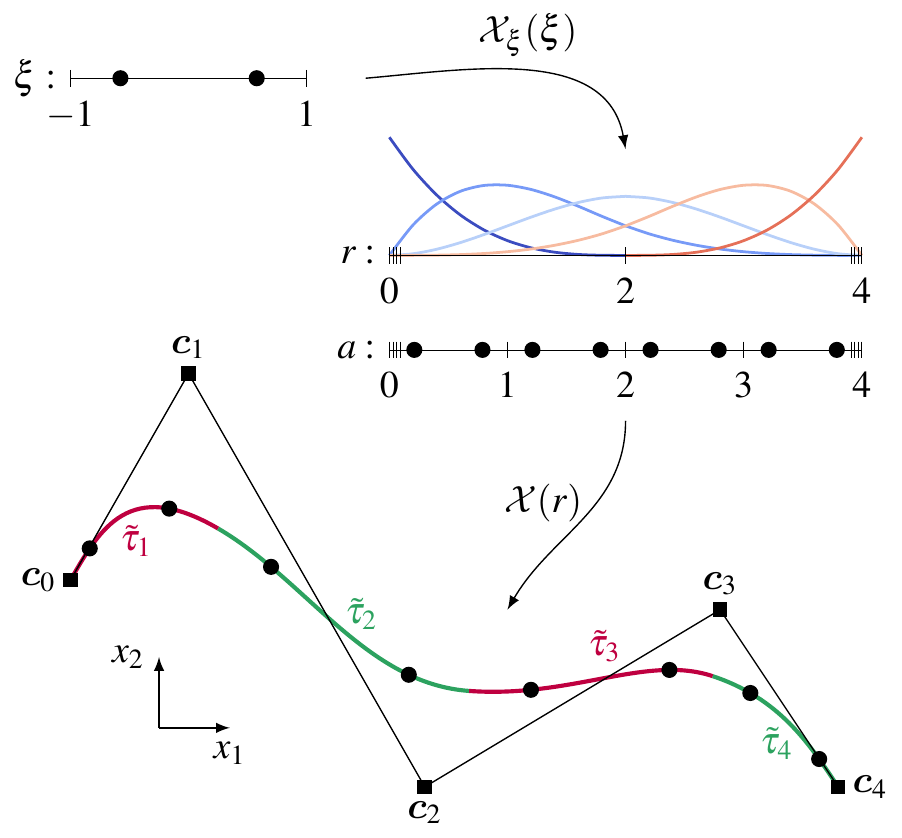}

In order to retain higher order convergence on domains with mixed boundary conditions or corners, the grading of elements towards such geometric locations may become necessary \cite{atkinson1997}. Corners are easily identified in the knot vector $\KVC$ multiplicity of knots $\multi=\pu$. The grading is performed by subdivision of the adjacent knot spans $\indexSpan=[\uuc_{\indexSpan},\uuc_{\indexSpan+1})$ 
into $m$ elements. This is simply performed in the parameter space by means of knot insertion in $\KVC$. The resulting $n$ knot values $\bar{\uuc}_\indexA$ with $\indexA = \fromto{1}{n}$ are defined by
\begin{align}
  \label{eq:MeshGrading}
  \bar{\uuc}_\indexA &= \uuc_{\indexSpan}  + \left(\uuc_{\indexSpan+1} - \uuc_{\indexSpan}\right)%
  \left( \frac{ \indexA}{n} \right)^{q} && \und && %
  \bar{\uuc}_\indexA = \uuc_{\indexSpan+1}  - \left(\uuc_{\indexSpan+1} - \uuc_{\indexSpan}\right)%
  \left( \frac{ \indexA}{n} \right)^{q}
\end{align}
for grading towards $\uuc_{\indexSpan}$ or $\uuc_{\indexSpan+1}$ respectively. 
The exponent $q$ is defined by
\begin{align}
  q &\geq \frac{\puc+1}{\hoelder} && \with && 0 < \hoelder \leq 1
\end{align}
where $\puc$ is the order of the quadrature rule and $\hoelder$ denotes the Hölder constant \cite{atkinson1997}.

In the current implementation, mixed boundary conditions are considered insofar, that either Dirichlet or Neumann boundary conditions along a single patch are allowed. Grading is performed in the vicinity of patches with different boundary conditions. 

\subsubsection{Local Refinement for Tensor Product Surfaces}
\label{sec:iga:refinement}

The use of artificial knot vectors $\KVC$ is sufficient for the partition of a curve into elements $\ieg$. However, the extension of this concept to surface representations is limited. In particular, local refinement is not possible if elements $\ieg$ are defined by a tensor product of $\KVC_\totalA$ and $\KVC_\totalB$. In this section, a strategy for local refinement with the isogeometric Nyström method is explained. The procedure is visualized in \figref{fig:refinementpoints}. 

Global refinement is performed by means of knot insertion, i.e. inserting $\bar{\uuc_\indexA}$ in $\KVC_\totalB$ as indicated by the dashed line in \mbox{\figref{fig:refinementpoints}(a)}. Both non-zero knot spans in $\KVC_\totalA$ are subdivided. A further subdivision is indented to be local for each of the elements $\ieg$. For the partition of an element $\ieg$ into \emph{local elements} $\iegLocal$ the definition of \emph{refinement points} $\locRefInfo$ is adequate. %
Each $\locRefInfo$ is defined in the parametric space and located either inside or on the edge of $\ieg$. Inside $\ieg$ the refinement point defines the origin of a cross which is aligned to the parametric coordinate system and subdivides $\ieg$ into four local elements $\iegLocal$. If $\locRefInfo$ is located on the edge, $\ieg$ is subdivided into two $\iegLocal$. Further, a local grid can be defined by combining several refinement points simultaneously.
The described local refinement options are illustrated in \mbox{\figref{fig:refinementpoints}(b)}.

In order to enable further refinement of local elements as well, all local elements are sorted in a \emph{hierarchical} tree structure and labeled with the refinement level $\indexLevel$. The initial refinement level $\indexLevel=0$ refers to the global element, hence $\iegLocal^0=\ieg$. Each node of the tree may have a different number $I$ of ancestors because of the manifold possibilities in defining $\locRefInfo$ per level. The local elements $\iegLocal^{\ell}_{\indexA}$ generated in level $\ell$ cover the complete area of the local element $\iegLocal^{\ell-1}$ of the previous level
\begin{align}
  \label{eq:NestedLocalElement}
  \iegLocal^{\ell-1} & = \bigcup_{\indexA}^{\totalA} \iegLocal^{\ell}_{\indexA}.
\end{align}
\tikzfig{tikz/refinementpoints}{2.2}{Subsequent refinement of elements defined by the tensor product of $\KVC_\totalA = \left\{0,0,0,1,1,2,2,2\right\}$ and $\KVC_\totalB = \left\{0,0,0,2,2,2\right\}$. (a) Global refinement by inserting $\bar{\uuc_\indexA} = 1$ into $\KVC_\totalB$. (b) Subdivision into local element by refinement points of the first refinement level which are denoted by circles. (c) Further local refinement by a higher refinement level indicated by diamonds.
}{fig:refinementpoints}{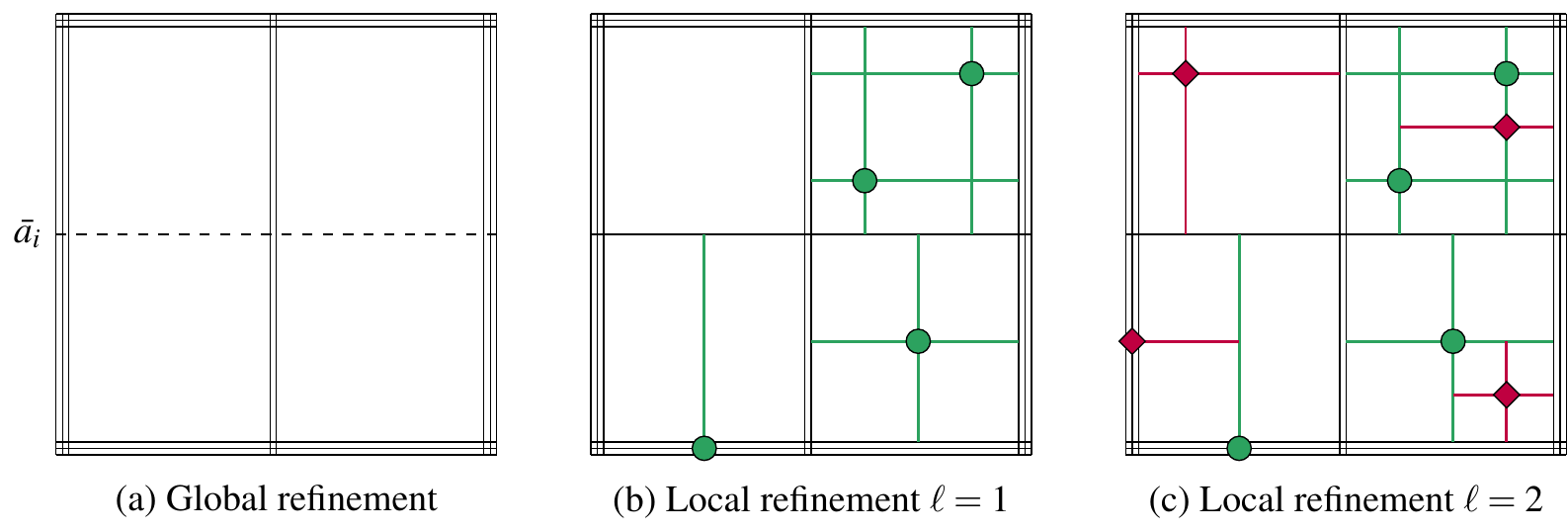}

The final partition of the global element $\ieg$ is defined by the sum of $\totalB$ local elements
\begin{align}
  \label{eq:partitionOfLocalElement}
  \ieg &= \bigcup_{\indexB=1}^{\totalB} \iegLocal_{\indexB}
\end{align}
related to the \emph{leafs} in the hierarchical tree structure. An example of such a locally refined patch with two levels of refinement is depicted in \mbox{\figref{fig:refinementpoints}(c)}. The local refinement procedure involves the scaling and translation of the element boundaries. Details on the construction of this mapping due to a given set of $\locRefInfo$ are found in \ref{appsec:appendixTransformation}.

\subsection{Isogeometric Nyström Method with Local Correction}
\label{sec:iga:nystroem}

In the presented implementation, Gauss-Legendre quadrature rules are taken. For the analysis in three dimensions, a tensor product quadrature is constructed as illustrated in \figref{fig:LocalRefinement}. However, it is also feasible to apply non-tensor product quadrature or numerical quadrature constructed for special special purposes in that context~\cite{bremer2010c}. 

As mentioned in \secref{sec:bie:bie}, the integral kernels are singular if collocation and quadrature point coincide $\pt{x}_i=\pt{y}_j$. Such kernels require special treatment for a correct integration. For the isogeometric Nyström method, a spatial separation of quadrature points in relation to each collocation point is performed. 
An admissibility criterion
\begin{align}
  \label{eq:AdmissibilityCriterion}
  \diam \left( \iegLocal \right) & \leq \eta \dist \left( \pt{x}_i, \pt{y}_j \right)
\end{align}
is introduced which separates the regime with a smooth kernel function from that one with singular or nearly singular behavior. 
If \eqref{eq:AdmissibilityCriterion} is fulfilled, the corresponding entries are in the scope of the \emph{far field} where the system matrices consist of point evaluations only. In particular, a matrix entry of the discrete single layer potential \eqref{eq:SL} for the isogeometric Nyström method reads
\begin{align}
  \label{eq:SL_IGA}
  \mat{V}[i,j] = \fund{U}(\pt{x_i},\pt{y}_j) G(\pt{\uu}_j) J_{\mychi_\xi}(\pt{\xi}_j)  w_j. 
\end{align}
In \eqref{eq:SL_IGA}, $\fund{U}(\pt{x_i},\pt{y}_j)$ is the evaluation of the fundamental solution with respect to the spatial coordinate of the collocation point $\pt{x_i}$ and that of the quadrature point $\pt{y}_j$ in $\R^d$.
%
\tikzfig{tikz/quadratureDistributionLocalElement}{1.0}{The distribution of the quadrature points $\GP$ on a locally refined surface. The basis function of the geometry are defined by $\KV_\totalA = \KV_\totalB = \left\{0,0,0,2,2,2\right\}$. The partition of local elements is defined by a global insertion of $\bar{\uuc}_1 = \bar{\uuc}_2 = 1$ and a refinement point $\locRefInfo = (0.5,0.5)^\trans$.}{fig:LocalRefinement}{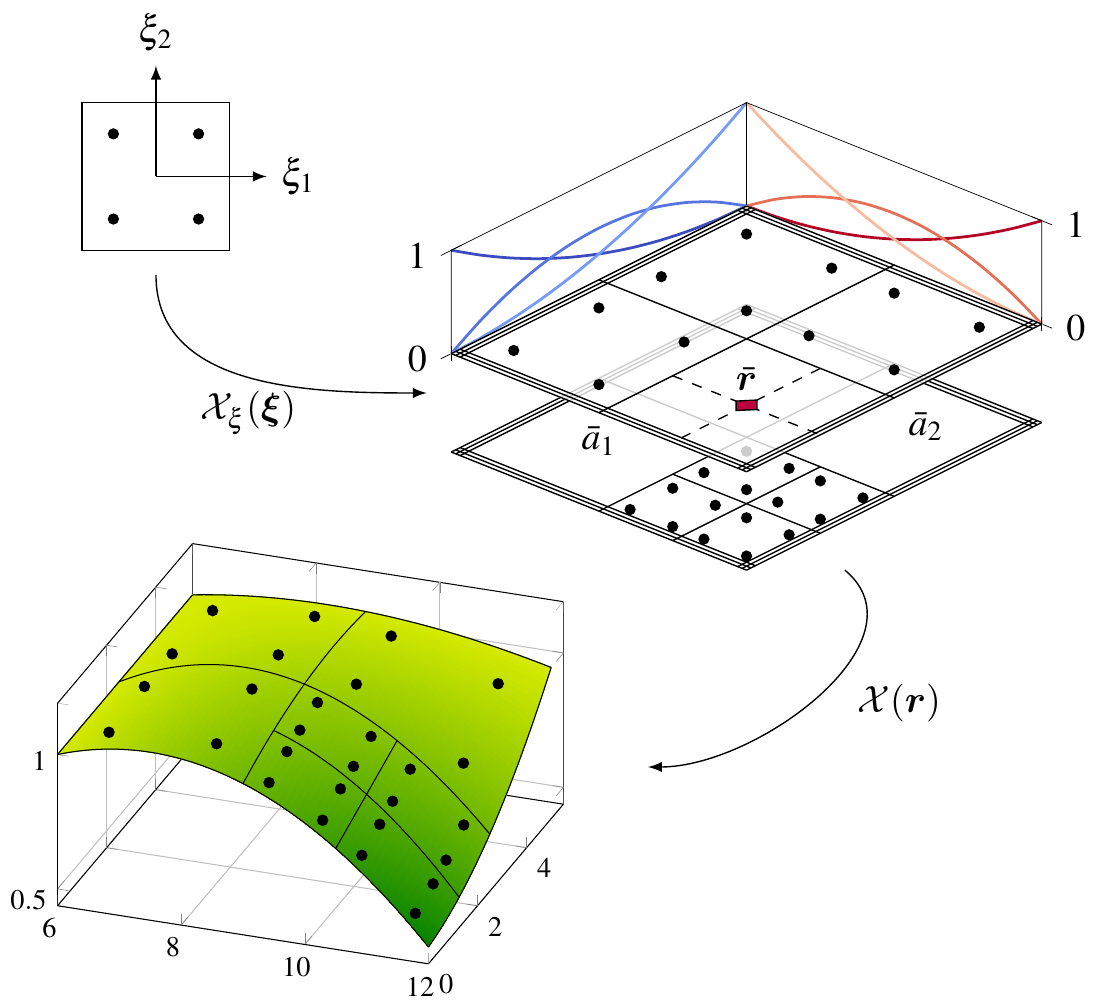}
The evaluation of Gram's determinant 
\begin{equation}
  \label{eq:GramsDeterminant}
  G(\pt{\uu}_j)=\sqrt{\det\left(\mat{G}(\pt{\uu}_j)\right)}
\end{equation}
is defined by the location $\pt{\uu}_j$ of the quadrature point $\pt{y}_j$ in the parametric space. Since the quadrature points are given in the reference element by its coordinates $\pt{\xi}_j=(\xi_{j,1},\xi_{j,2})^{\trans}$, the mapping $\pt{\uu}_j=\mychi_{\xi}(\pt{\xi}_j)$ to the local element $\iegLocal$ is sufficient. For the integral transformation from reference to local element, the Jacobian of the mapping $\mychi_\xi$
\begin{equation}
  \label{eq:Jacobian}
  J(\pt{\xi}_j) = \det\left( \mat{J}(\pt{\xi}_j) \right)
\end{equation}
is evaluated with respect to the reference coordinates $\pt{\xi}_i$ of the quadrature point $\pt{y}_i$. Finally, $w_j$ in \eqref{eq:SL_IGA} denotes the original quadrature weight.

The \emph{near field} zone, where the integral kernels are singular or nearly singular along the affected local elements is defined for point combinations where \eqref{eq:AdmissibilityCriterion} is not fulfilled. At that point, the matrix evaluations are performed by means of the locally corrected integral kernels as described in \secref{sec:bie:nystroem}. To perform the local correction, B-splines are taken for the polynomial test functions $N_i = \Bspline_{i,\pu}$ in \eqref{eq:LocalCorrection} and defined on the local element $\iegLocal$. In particular, a \BEZ interpolation is proposed which is defined by the knot vectors of the kind
\begin{align}
	\label{eq:knot_vector_BEZ_local_correction}
	\KV &= \left\{ {-1}, \dots, {-1}, 1, \dots, 1 \right\} 
\end{align}
in all parametric directions. The multiplicity of the knots is chosen so that they define at least as many basis functions as present quadrature points on the local element. This allows the solution of \eqref{eq:LocalCorrection} by means of $LU$-decomposition or by solving a least squares problem.

The integrals of the right hand side in \eqref{eq:LocalCorrection} are
\begin{align}
  \label{eq:SL_LCN}
  \int_{\iegLocal} \fund{U}(\pt{x}_i,\pt{y}) \Bspline_{i,\pu}\ofpt{y} \dgamma{y} && \und && %
  \int_{\iegLocal} \dTr_{\pt{y}} \fund{U}(\pt{x}_i,\pt{y}) \Bspline_{i,\pu}\ofpt{y} \dgamma{y}.
\end{align}
If $\pt{x}_i \in \iegLocal$ the integrals are weakly or strongly singular. In that case the single layer integral is subject to transformation described in \cite{lachat1976} and \cite{duffy1982} while the double layer integral is treated with regularization techniques presented in \cite{guiggiani1990} respectively. If \eqref{eq:AdmissibilityCriterion} is not fulfilled but $\pt{x}_i \not\in \iegLocal$, then the integral is nearly singular and treated with adaptive numerical integration as described in \cite{marussig2015}. Practically, the extent of the region where local elements are marked as nearly singular is determined by the admissibility factor $\eta$.



\subsection{Isogeometric Postprocessing}
\label{sec:iga:postprocessing}

Once the system of equations is solved, the Cauchy data exists only in the quadrature points $\GP$. Thus, post-processing steps are required to visualize the distribution over the whole geometry. The Nyström-interpolation \cite{atkinson1997} is the most accurate procedure for this task. But it requires additional kernel evaluations at all quadrature points, which is computationally expensive. For the isogeometric Nyström method, the following approach is probably less accurate but simpler and local to $\iegLocal$. Following the isogeometric concept, each element $\iegLocal$ is represented by the \BEZ patch already constructed for the local correction. The results in each quadrature point $\pt{y}_j$ are interpolated within each $\iegLocal$ by means of the basis functions $\Bspline_{\indexA,\pu}$ based on a knot vector of form~\eqref{eq:knot_vector_BEZ_local_correction}. For instance, the primary variable in any point $\pt{\xi}$ on the reference element can be calculated with 
\begin{align}
	\label{eq:resultsInterpolation}
	u\ofpt{\xi} &= \sum_{\indexA=0}^{\totalA} \Bspline_{\indexA,\pu}\ofpt{\xi} \: \CP_\indexA.  
\end{align}
In order to compute the unknown coefficients $\CP_\indexA$ the inverse of the mapping $\mat{C}\ofpt{\xi}\vek{c} = \vek{u}$ is needed which is defined by the \emph{spline collocation matrix}. Its entries are
\begin{align}
  \label{eq:splineCollocationMatrix}
  \mathbf{C} [\indexB,\indexA] &=  \Bspline_{\indexA,\pu}(\pt{\xi}_\indexB) && \mit && %
  \indexA=\fromto{0}{\totalA} && \und && \indexB=\fromto{0}{\totalB}
\end{align}
where $\indexA$ is the number of B-spline functions and $\indexB$ the number of quadrature points on $\iegLocal$. The linear system~\eqref{eq:resultsInterpolation} can be solved directly or in a least squares sense.


%% file: results.tex
\section{Numerical Results}
\label{sec:results}

In this section, numerical results are provided for academic and practical problems. The results are critically reviewed and remarks on limitations and open topics of the isogeometric Nyström method are given.


For the numerical analysis of the convergence of the isogeometric Nyström method, problems with an infinite domain $\Omega$ are solved. 
The fundamental solution $\fund{U}(\tilde{\pt{x}},\pt{y})$ with a number of source points $\tilde{\pt{x}}\in\Omega^{-}$ outside of $\Omega$ is applied as a boundary condition at the quadrature points $\pt{y}\in\Gamma$. 
The problem is solved by means of Fredholm integral equations of the first kind \eqref{eq:FredholmFirst} in order to test the single layer potential $\mat{V}$ and with the second kind equation \eqref{eq:FredholmSecond} to test the double layer potential $\mat{K}$.
Results are given in the interior at several points $\hat{\pt{x}}$ and the error is defined by
\begin{align}
  \label{eq:indirect-check}
  \err_h & = \primary(\hat{\pt{x}}) - \fund{U}(\tilde{\pt{x}},\hat{\pt{x}})  && \forall\hat{\pt{x}}\in\Omega \komma \tilde{\pt{x}}\in\Omega^{-}.
\end{align}
The relative error is 
\begin{equation}
  \label{eq:relative_error}
  \err_{rel} = \frac{\err_h}{\fund{U}(\tilde{\pt{x}},\hat{\pt{x}})}
\end{equation}
and measured in the maximum-norm $\|\err_{rel}\|_{\infty}$. The normalized element diameter
\begin{equation}
  \label{eq:meshwidth}
  h = \left( \frac{ A^{\ieg}_{max}}{A} \right)^{1/(d-1)}
\end{equation}
is used for convergence plots where $A^{\ieg}_{max}$ is the largest length or area of all elements $\ieg$ and $A$ the surface length or area of the whole boundary $\Gamma$. Hence, a step of the process denoted as \emph{uniform $h$-refinement} halves the knot span of $\ieg$ in two dimensions ($d=2$) resulting in two new elements. In three dimensions ($d=3$), the knot spans in both parametric directions are affected which produces four new elements. However, the parameter~$h$ always refers to the length or area in $\R^2$ or $\R^3$ respectively. 

In the following sections, the terminus \emph{$p$-refinement} refers to the step-wise increase of the quadrature order used for the simulation. For one-dimensional elements $\ieg$ in $\R^2$, a step of $p$-refinement results in the increase by one of the local quadrature points per element. For the tensor product quadrature being used for analysis in $\R^3$, this process leads to $2p-1$ additional quadrature points on $\ieg$.

  \input{\CommonPath/results_flower}
\clearpage

  \input{\CommonPath/results_teardrop}
\clearpage

  \input{\CommonPath/results_cantilever}
\clearpage

  \input{\CommonPath/results_torus}
\clearpage

  \input{\CommonPath/results_fichera}
\clearpage



%
  \input{\CommonPath/results_spanner}
\clearpage


%% file: results_flower.tex
\subsection{Flower}%
\label{sec:res:flower}
%
The convergence of the method in two dimensions is tested on a smooth flower-like geometry solving an exterior Dirichlet problem. The whole setting is depicted in \figref{fig:geometries2d} on the left. The admissibility factor for the local correction is set to $\eta=2.0$.%
\tikzfig{tikz/geometries2d.tex}{0.8}{Geometries and the corresponding NURBS control grid (squares) for the test settings in two dimensions: flower-like boundary in infinite domain with source points $\tilde{\pt{x}}$ and test points $\hat{\pt{x}}$ (left) and teardrop with opening angle $\alpha$ (right)}{fig:geometries2d}{figure07}%
In \figref{fig:flower_convergence} the results for Laplace's equation with uniform $h$-refinement are shown for different orders $p$ of the chosen quadrature rule. While the double layer operator $\mat{K}$ shows the expected optimal convergence of $p+1$, the single layer $\mat{V}$ does not. In that case, the convergence reverts to a linear rate irrespective of the order of the integration.%
\tikzfig{tikz/flower.xml_potential_pointwise}{1.0}{Convergence rates of the indirect Laplace problem for the discrete double layer (left) and single layer potential (right) on the two dimensional flower like geometry. The optimal convergence rates for the lowest and highest order are indicated by triangles.}{fig:flower_convergence}{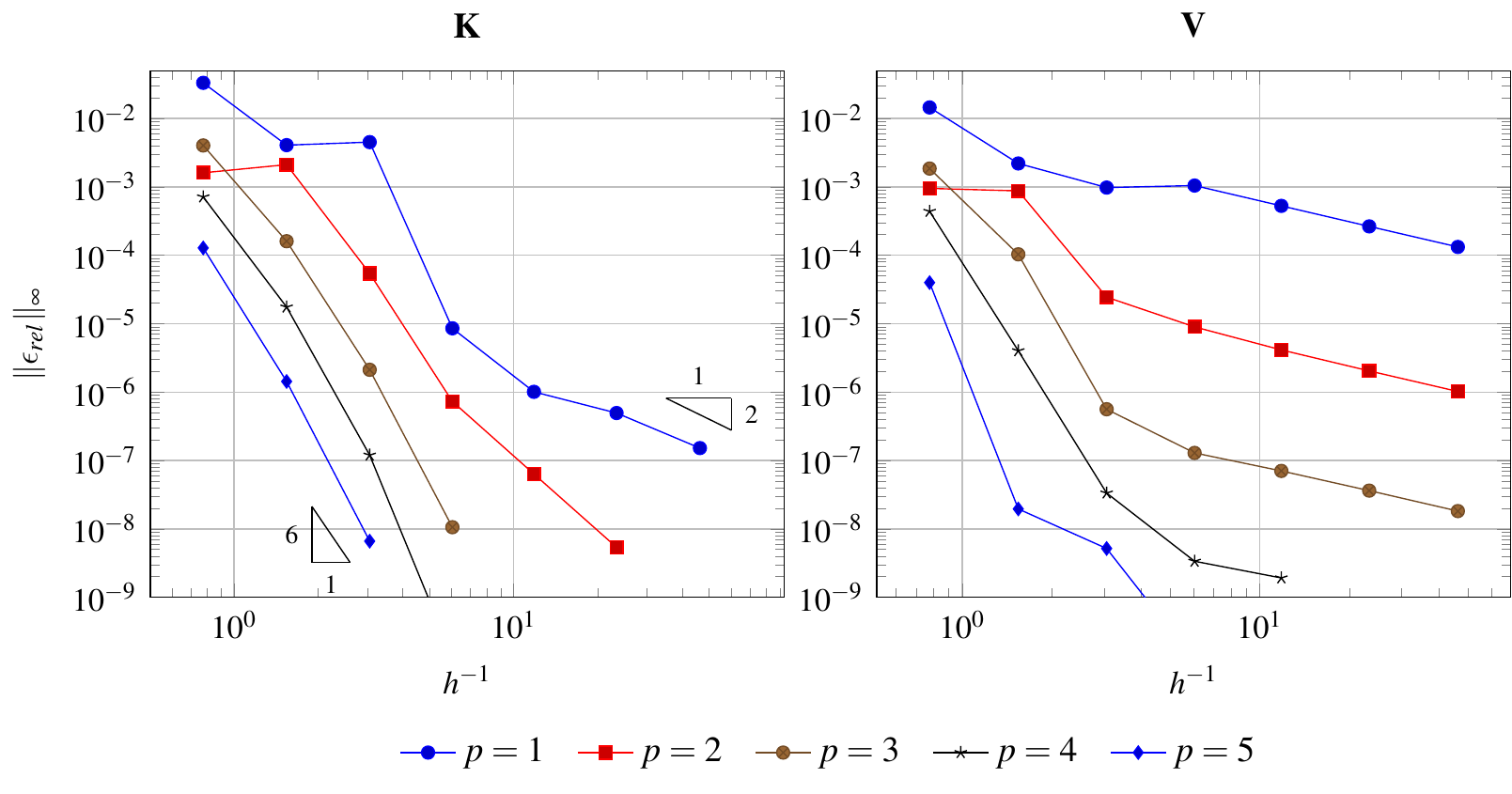}%
However, there is a significant offset depending on the quadrature rule so that a convergence following
\begin{equation}
  \label{eq:exp_convergence}
  \| \err_{rel} \|_{\infty} = C\exp^{(-n^{s})}
\end{equation}
 with respect to the number of degrees of freedom $n$ is observed for the $p$-refinement strategy. In the case of the flower-like geometry the factor $s$ is approximately $0.55$. The constant $C$ depends on the size of the integration elements. This is shown in \figref{fig:flower_convergence_exponential}, where $h_s$ denotes the initial, unrefined and dimensionless element size in the sense of equation~\eqref{eq:meshwidth}. The dashed lines follow~(\eqref{eq:exp_convergence}) with different exponential factors~$s$%
\tikzfig{tikz/flower.xml_potential_pointwise_exponential.tex}{1.0}{Convergence rates for $p$-refinement of the indirect Laplace problem for the discrete double layer (left) and single layer potential (right) on the two dimensional flower like geometry. The dotted lines indicate the exponential function~\eqref{eq:exp_convergence} where the aligned number is the exponential factor $s$}{fig:flower_convergence_exponential}{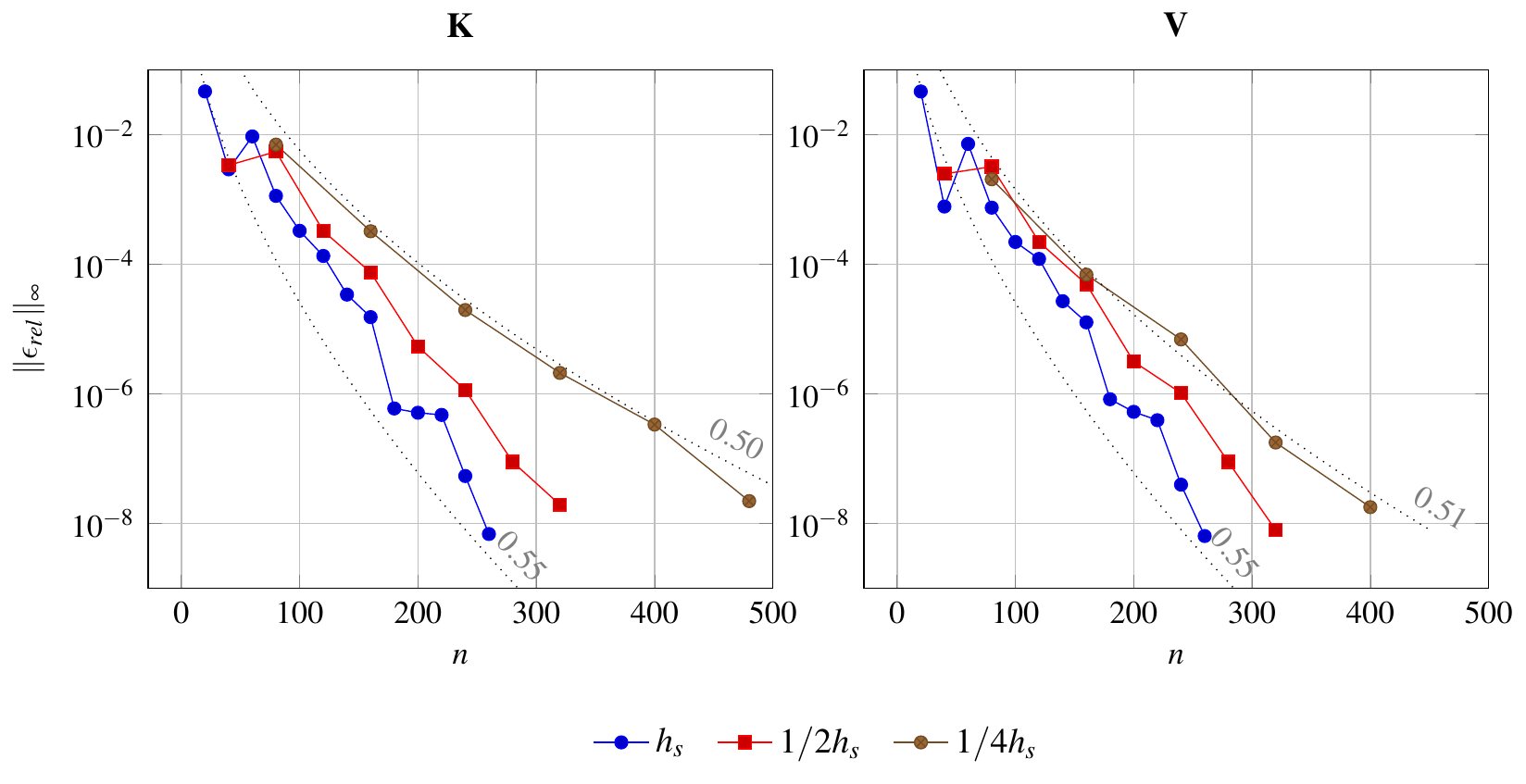}
For the \Lame-Navier equation, the behavior is depicted in \figref{fig:flower_convergence_lame} where the \Lame-parameter are defined by the elastic constants
\begin{align}
  \label{eq:elasticConstants}
  \lambda=\frac{E\nu}{(1-2\nu)(1+\nu)} & & \und & &
  \mu=\frac{E}{2(1+\nu)}.
\end{align}
In that example, the Young's modulus and Poisson's ratio are set to $E=\SI{1}{\giga\pascal}$ and $\nu=0.3$ respectively.%
\tikzfig{tikz/flower.xml_planestress_pointwise}{1.0}{Convergence rates of the indirect \Lame-Navier problem for the discrete double layer (left) and single layer potential (right) on the two dimensional flower like geometry.}{fig:flower_convergence_lame}{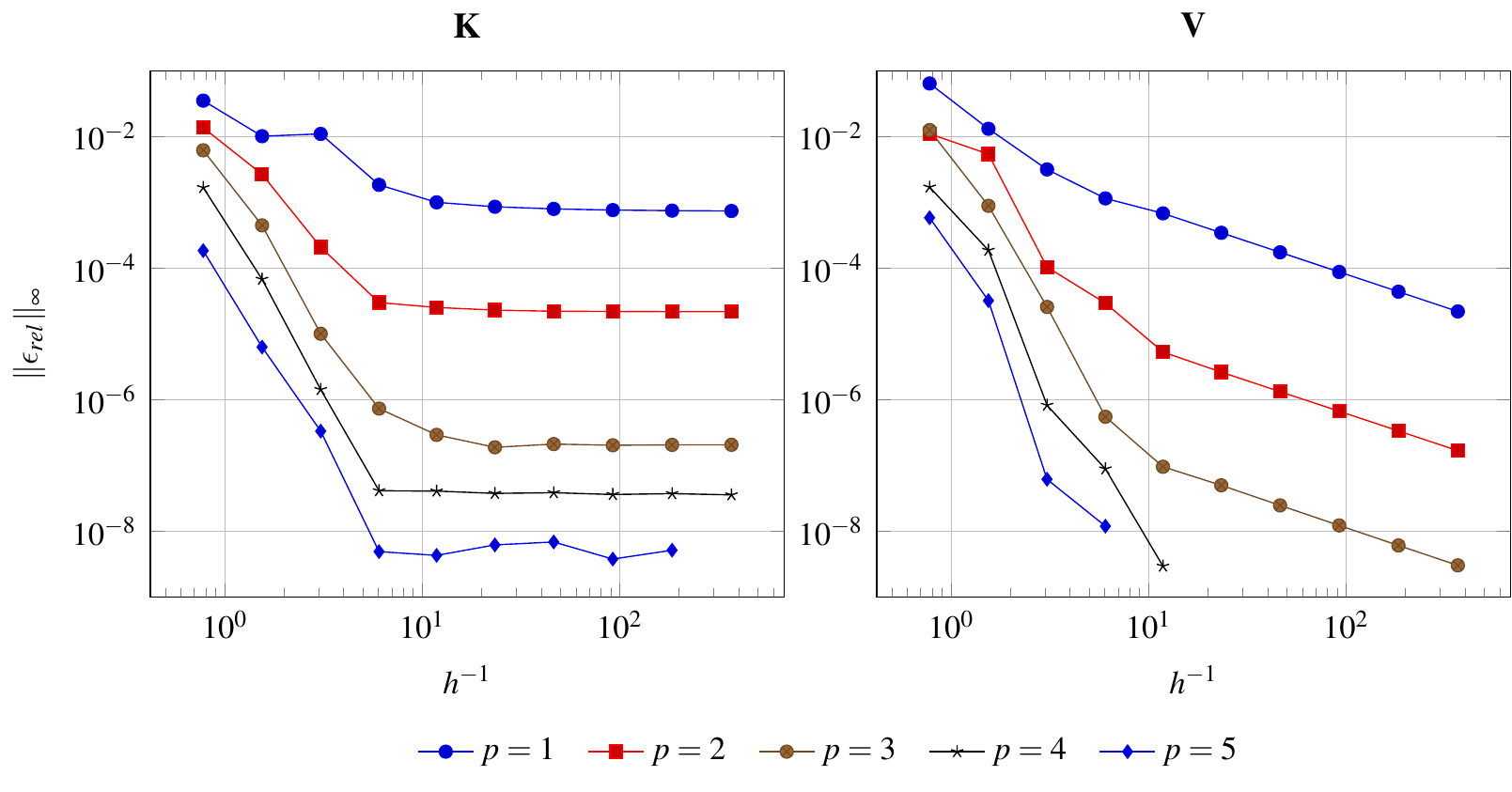}
Additionally, \figref{fig:flower_convergence_lame_exponential} depicts the exponential convergence with respect to the degrees of freedom $n$ for the $p$-refinement strategy, in that case with $s\approx 0.50$.
\tikzfig{tikz/flower.xml_planestress_pointwise_exponential}{1.0}{Convergence rates for $p$-refinement of the indirect \Lame-Navier problem for the discrete double layer (left) and single layer potential (right) on the two dimensional flower like geometry. The dotted lines indicate the exponential function~\eqref{eq:exp_convergence} where the aligned number is the exponential factor $s$}{fig:flower_convergence_lame_exponential}{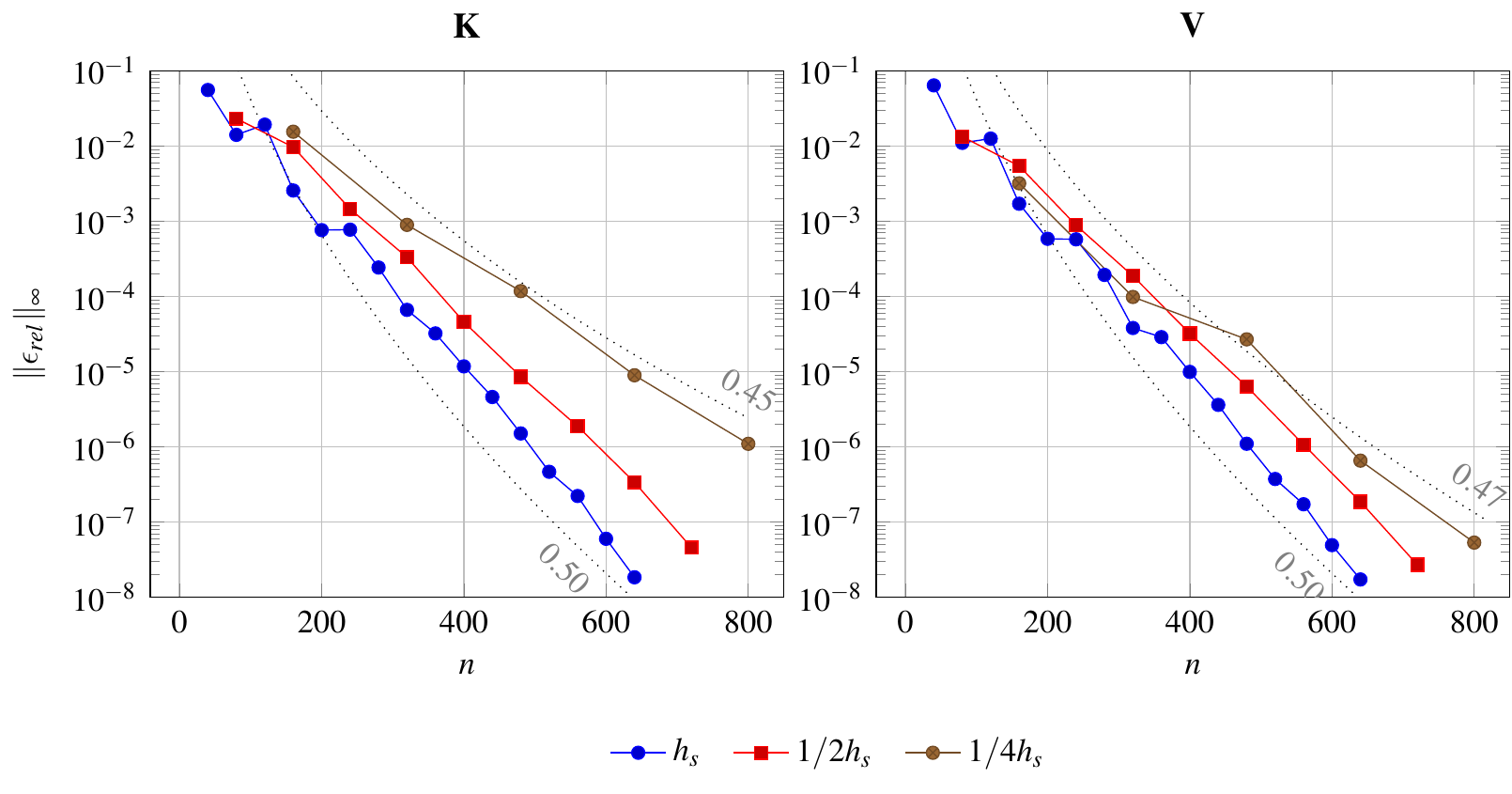}

For Fredholm integral equations of the second type~\eqref{eq:FredholmSecond} solved with the Nyström method, convergence theorems have been proven for the Laplace equation in two and three dimensions \cite{atkinson1997}. The numerical results for the double layer operator in \figref{fig:flower_convergence} comply with the theory and show convergence rates of $p+1$ or better for a quadrature rule of order $p$.
However, for boundary integral formulations in linear elasticity with the \Lame-Navier equation, literature is rather sparse. In the straight-forward implementation described in this paper, convergence for the discrete double layer operator stops at a certain level of $h$-refinement.

In contrary to the Laplace equation, where the double layer operator is weakly singular, the boundary integral for problems in elasticity is evaluated in the sense of a Cauchy principal value which requires particular attention. In the presented implementation it is ensured, that the accuracy for the local correction is beyond the measured relative error $\|\err_{rel}\|_{\infty}$. This affects the right hand side in equation~\eqref{eq:LocalCorrection}. 
Numerical studies have shown, that a reasonable variation of the admissibility factor $\eta$ does not improve the convergence significantly. An additional study, which has not been presented in this paper for brevity, indicated that if the admissibility factor is chosen to be $\eta\rightarrow\infty$ and hence every integration element is corrected, the numerical results show optimal convergence agin. But this also means that the isogeometric Nyström method is in fact replaced by the isogeometric BEM for which convergence has been demonstrated numerically in \cite{marussig2015}. However, it is remarkable that the convergence plateau for the flower example is reached after exactly $3$ uniform refinement steps for all quadrature orders as depicted in \figref{fig:flower_convergence_lame}. By altering the element partition such that $h_{max}\approx h_{min}$ the kink has been shifted to one more refinement step, but still the plateau does not disappear. 

The boundary integral operator for Fredholm integral equations of the first kind is not compact anymore and hence, mathematical analysis is rather involved  \cite{sauter2011}. Moreover, for problems in two dimensions the fundamental solution is logarithmic. The direct application of the presented framework to the boundary integral equation results in only linear convergence independent of the chosen quadrature. The performed numerical tests show this behavior for both the Laplace and \Lame-Navier problem. For logarithmic first kind boundary integral equations on closed curves in $\R^2$ it is proposed in literature to apply a transformation so that they become harmonic. Therefor, convergence can be restored with respect to the numerical quadrature \cite{sloan1981}.

Although the convergence with $h$-refinement for two dimensional problems is restricted, there is a significant offset between discretizations with increasing quadrature order. As a consequence, the $p$-refinement strategy is preferred over $h$-refinement for the analysis of practical applications with the isogeometric Nyström method. Convergence follows equation~\eqref{eq:exp_convergence} in that case.

%% file: results_teardrop.tex
\subsection{Teardrop}
\label{sec:res:teardrop}

For the discretization of the direct BIE~\eqref{eq:DirectBIE} with the Nyström method, singularities of the solution occur at domains with corners. This is, because the gradient of the solution becomes singular. For the isogeometric Nyström method the problem is tackled by the grading strategy described 
in \secref{sec:iga:cauchy}. 

To validate the algorithm, a two dimensional domain modeled like a teardrop is analyzed. The domain, which is shown on the right of \figref{fig:geometries2d}, has a single corner where the opening angle is chosen to be $\alpha=\SI{90}{\degree}$. The second kind integral equation~\eqref{eq:IndirectMatrixDLP} for the exterior Laplace problem is solved and it's convergence behavior investigated. In a first step, uniform $h$-refinement is performed for different quadrature orders. Then, grading towards the corner is applied which is specified by equation~\eqref{eq:MeshGrading} where the number of sub-elements is chosen to be $m=6$ and the Hölder constant is $\hoelder=1$. To handle the small elements near the corner, the admissibility factor for the local correction is set to $\eta=6.0$. Both results are depicted in \figref{fig:teardrop_DLP}. The convergence is calculated with respect to degrees of freedom $n$.\mycomment[JZ]{Konvergenz ohne Grading: 2 Rechnungen linear fehlen}%
\tikzfig{tikz/teardrop_DLP_compare.tex}{1.0}{Convergence of the indirect Laplace problem on an infinite domain with one corner. The geometry is subject to uniform (left) and graded $h$-refinement (right). Convergence rates are indicated by triangles.}{fig:teardrop_DLP}{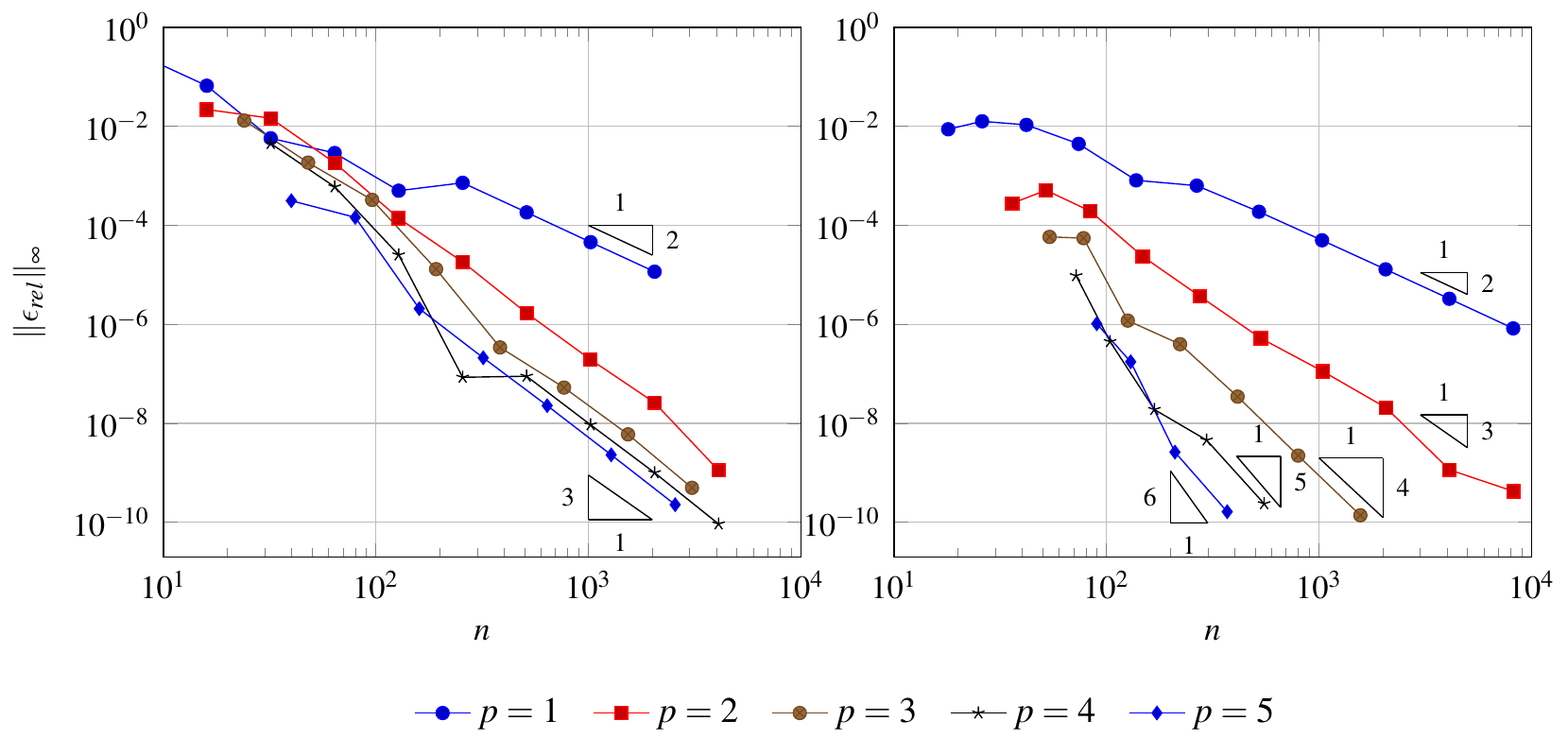}
%

The tests on the teardrop geometry justify the grading approach. In the case of the teardrop, the full convergence for Fredholm integral equations of the second kind is restored. As outlined in the introduction, several other techniques are available to tackle such singularities. All the mentioned approaches may be taken into account by the presented framework for the isogeometric Nyström method without any restriction.


%% file: results_cantilever.tex
\subsubsection{Cantilever}
\label{sec:res:cantilever}

\mytodo[BM]{%
  \begin{itemize}
  \item [\done]  Lame
  \item [\done]  Direct formulation
  \item [\done]  Mixed BC, corners
  \item [\cancelled]Compare BEM?
  \item Sketch of setup and displacement sketch from Paraview
  \end{itemize}
}

The last problem tested in two dimensions is a cantilever beam. The beam's dimensions are $\SI{10}{\meter}\times\SI{1}{\meter}$ is subject to a constant vertical top loading $t_y=\SI{-1}{\mega\pascal}$. The elastic constants are $E=\SI{29000}{\mega\pascal}$ and $\nu=0.0$. The solution with the isogeometric Nyström method is compared with the analytic solution by Timoshenko. The vertical displacement $u_y$ is measured at the midpoint $(10,0)$ of the free end of the cantilever. Grading of the integration elements towards the corners is realized with a Hölder constant $\hoelder=1$ and $m=4$ sub-elements. The admissibility factor for the local correction is chosen to be $\eta=6.0$. The vertical displacement for different quadrature orders and different degrees of freedom $n$ are presented in \figref{fig:cantilever_results} showing excellent results for all orders $p>1$.%
\tikzfig{tikz/cantilever2d_nystroem_LCN(6)_graded_h1_i4.tex}{1.0}{Vertical displacement $u_y$ in~[\si{\meter}] for the cantilever beam at its free end calculated with the isogeometric Nyström method. The analysis is carried out with different quadrature orders $p$ and degrees of freedom $n$. The analytic solution is indicated by the dashed line.}{fig:cantilever_results}{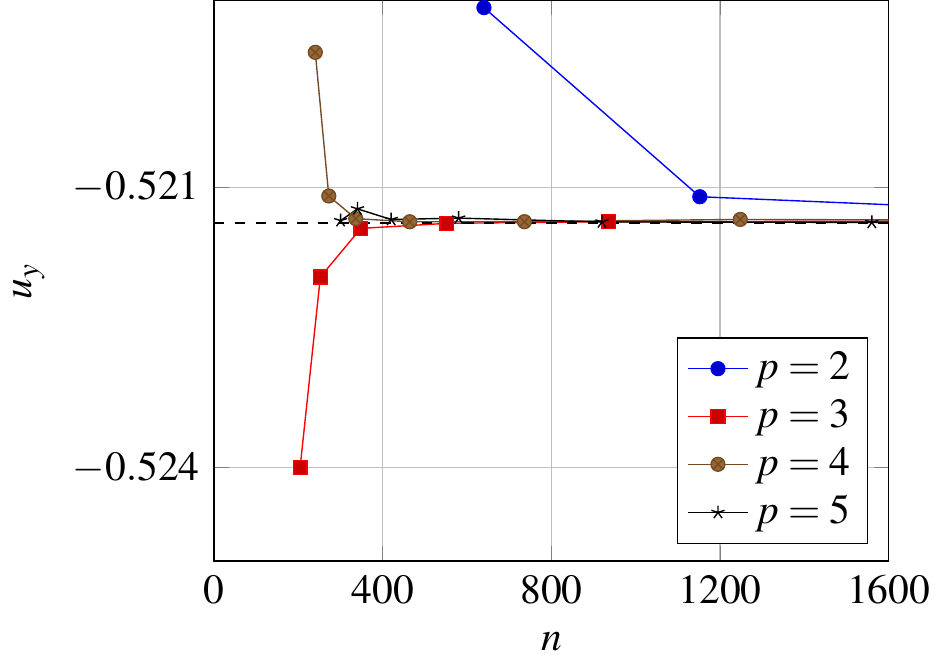}
%


%% file: results_torus.tex
\subsection{Torus}

\mytodo[JZ]{
  \begin{itemize}
  \item [\done] Laplace DLP h-refinement
  \item Laplace/Lame SLP h-refinement
  \item [\done] Lame DLP h-refinement
  \item Lame 5th-order calculation
  \end{itemize}
}
A torus geometry as depicted in \figref{fig:torus} is used to show the convergence of the boundary integral operators in three dimensions. The distance from the torus center to the center of the tube is $r_1=\SI{0.9}{\meter}$. The tube-radius is set to $r_2=\SI{0.2}{\meter}$. Both, the Laplace and the \Lame-Navier equation are solved. For the latter, the elastic constants are chosen to be $E=\SI{1.0}{\giga\pascal}$ and $\nu=0.3$. The results are outlined in \figref{fig:torus_laplace_convergence} and \figref{fig:torus_lame_convergence}.%
\tikzfig{tikz/torus.tex}{2.0}{Torus geometry taken for the investigation of convergence in three dimensions}{fig:torus}{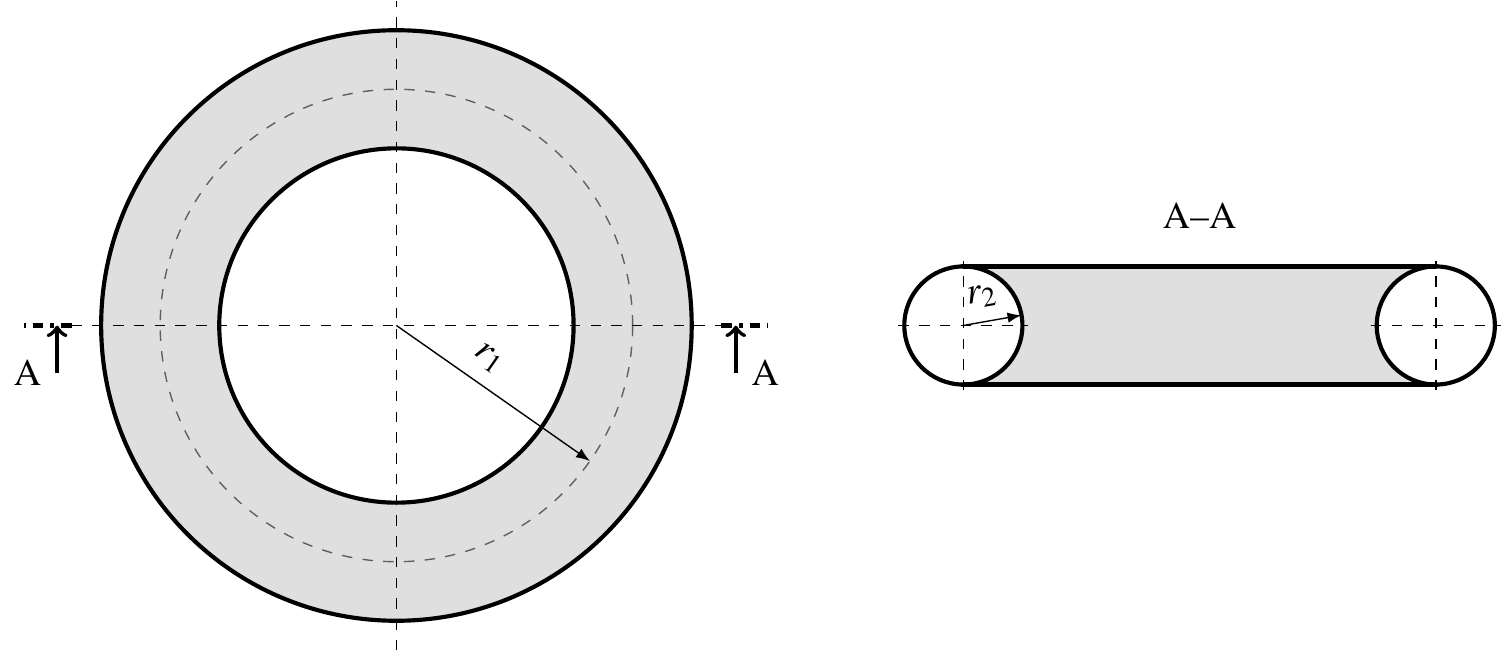}
\tikzfig{tikz/toruscustom.xml_potential_pointwise.tex}{1.0}{Convergence rates for the discrete double and single layer potential on the torus for the exterior Laplace problem.  The optimal convergence rates for the lowest and highest order are indicated by triangles.}{fig:torus_laplace_convergence}{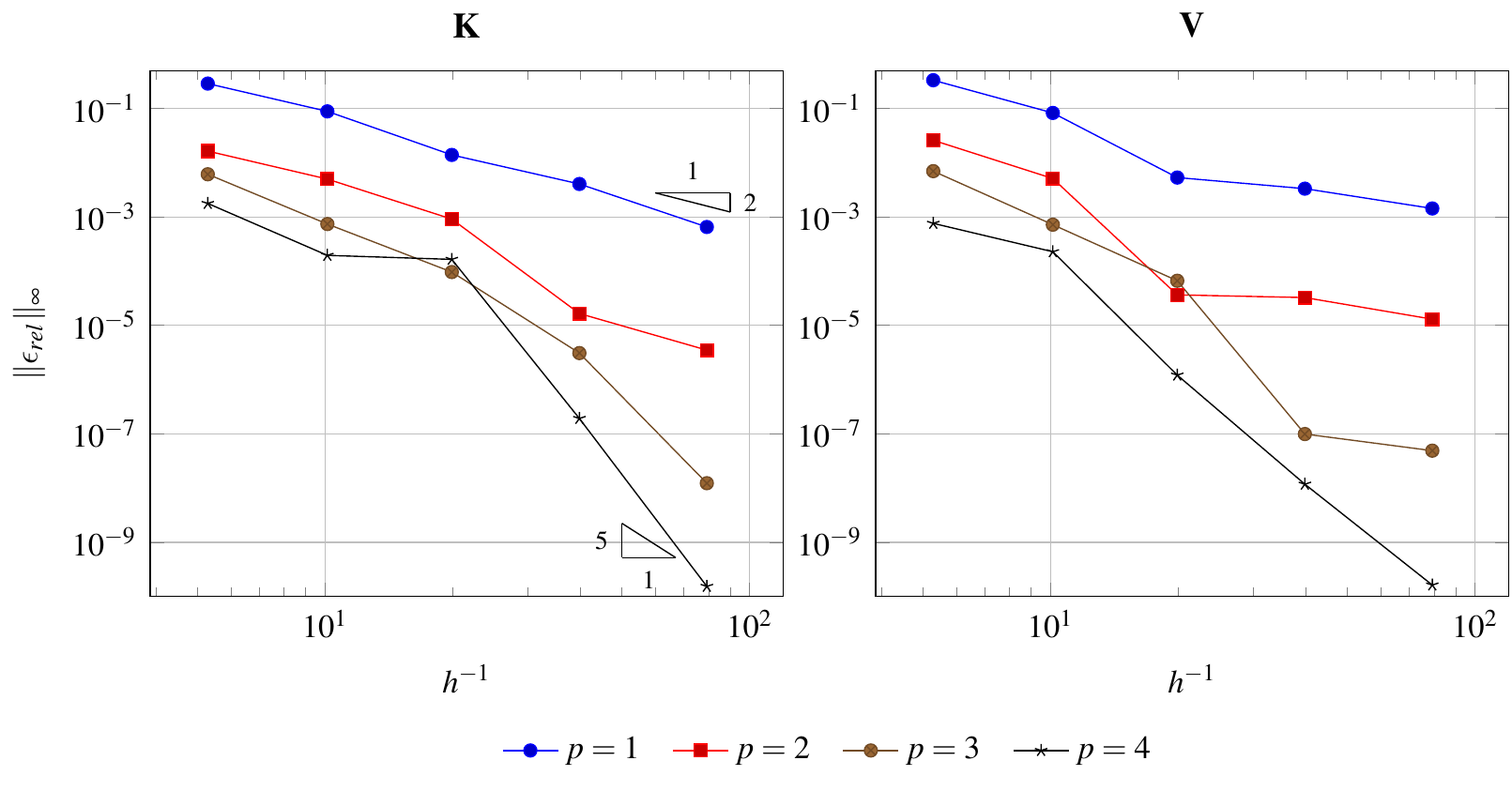}
\tikzfig{tikz/toruscustom.xml_lame_pointwise.tex}{1.0}{Convergence rates for the discrete double layer potential on the torus for the exterior elastic problem. Results are shown for the $h$-refinement (left) and the $p$-refinement strategy (right). On the left, the optimal convergence rates for the lowest and highest order are indicated by triangles. The dotted lines on the right indicate the exponential function~\eqref{eq:exp_convergence} where the aligned number denotes the exponential factor $s$}{fig:torus_lame_convergence}{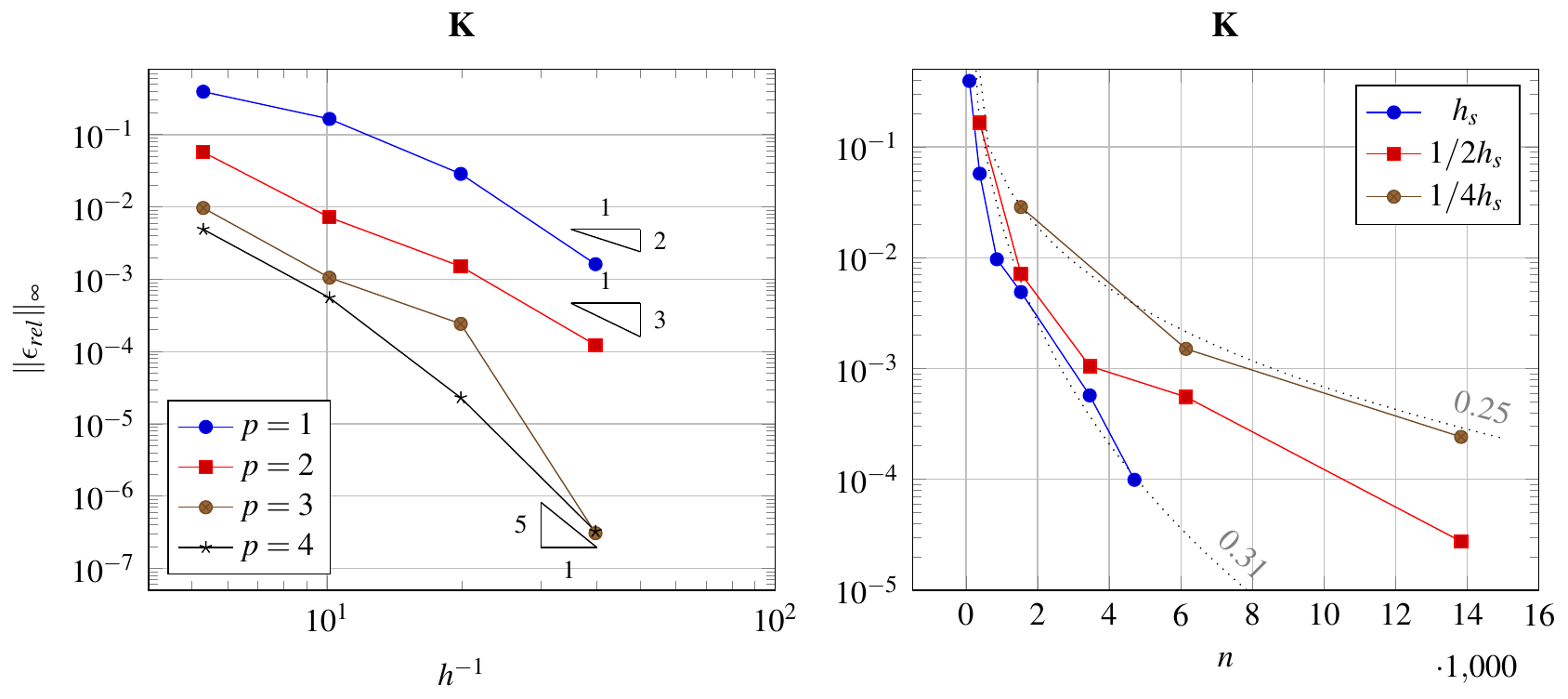}

For the Laplace equation, the numerical results in three dimensions comply with the theory as well. Similar to problems in two dimensions, optimal convergence $p+1$ for the double layer operator and suboptimal convergence for the single layer operator is observed. As for the convergence following equation~\eqref{eq:exp_convergence}, the factor is $s\approx 0.30$ for both, the discrete double layer operator $\mat{K}$ of the Laplace and the \Lame-Navier equation respectively. 

However, a convergence plateau such as for the flower geometry in \figref{fig:flower_convergence_lame} of \secref{sec:res:flower} has not been been observed. For further refinement steps the problem exceeded the number of degrees of freedom acceptable for such a study. The results for the three dimensional convergence study sketched in \figref{fig:torus_lame_convergence} required the application of a fast method which introduces additional numerical approximation. In the presented case, hierarchical matrices (\hmatrices) are taken, where the approximation quality of the boundary integral operators is controllable by a recursive criterion \cite{zechner2012c,zechner2013} and set to at least one magnitude lower than the achieved overall accuracy.


%% file: results_fichera.tex
\subsubsection{Fichera Cube}
\label{sec:res:fichera}

As a benchmark for the local refinement strategy, a heat conduction problem is solved on the Fichera cube and compared to an FEM solution \cite{rachowicz2006}.
The computational domain is given by $\Omega = (-1,1)^3 - [0,1]^3$ and represented by $14$ NURBS patches. On the Dirichlet boundary $\Gamma_D$ the temperature is set to ${g}_D = 0$. The flux on the Neumann boundary $\Gamma_N$ is an accumulation of known analytic solutions of L-shaped domains in two dimensions with respect to the $xy$-, $yz$- and the $xz$-plane. The temperature $u_{2d}$ and the heat flux $q_{2d}$ on the L-shape are given by
\begin{align}
\label{eq:l-shape-temp}
  \primary_{2d} &= r^{2/3} \sin\left( \frac{2\theta}{3} \right), && \with && r=\sqrt{\LshapeCoord_1^2+\LshapeCoord_2^2} && \und
\end{align}
\begin{align}
\label{eq:l-shape-flux}
  \dual_{2d} & =   \vek{n}[\indexA]  \frac{\partial \primary_{2d}}{ \partial \LshapeCoord_\indexA} , && \with && \indexA = 1,2
\end{align}
respectively. The local coordinates of the L-shape are denoted by $\LshapeCoord_\indexA$ and $\vek{n}$ represents the unit outward normal vector. The prescribed smooth flux ${g}_N$ on the Fichera cube is then accumulated from \eqref{eq:l-shape-flux} in each plane. The geometry of the problem is shown in \figref{fig:fichera_BC}.
\tikzfig{tikz/ficheracube.tex}{1.3}{Fichera cube geometry, where the Dirichlet boundary $\Gamma_D$ is drawn red and the Neumann boundary $\Gamma_N$ green. The right picture shows the  corresponding L-shape in two dimensions from which the non-zero Neumann boundary conditions are derived.}{fig:fichera_BC}{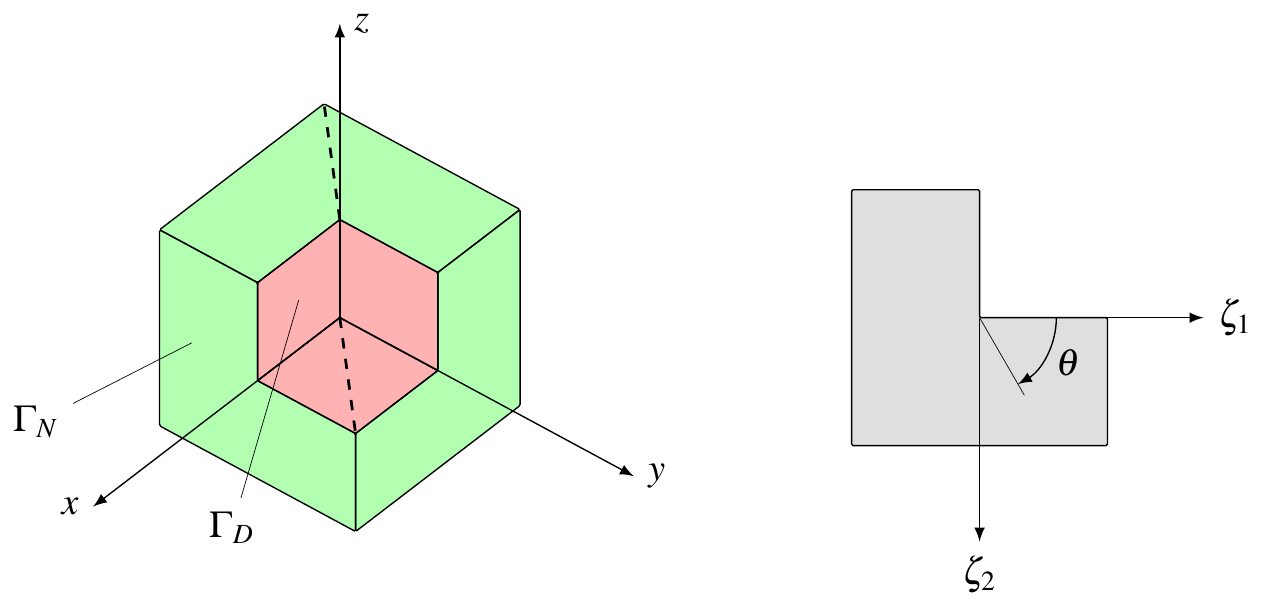}
The unknown Flux $\dual$ on $\Gamma_D$ has a singularity at the inverted corner of the Fichera cube. In order to approximate $\dual$ accurately towards this singular point, local refinement of integration elements as described in \secref{sec:iga:refinement} is applied. The first three refinement steps are illustrated in \figref{fig:fichera_grading}.%
\tikzfig{tikz/ficheracube_refinement.tex}{1.3}{The first three step of the local $h$-refinement on the Fichera cube}{fig:fichera_grading}{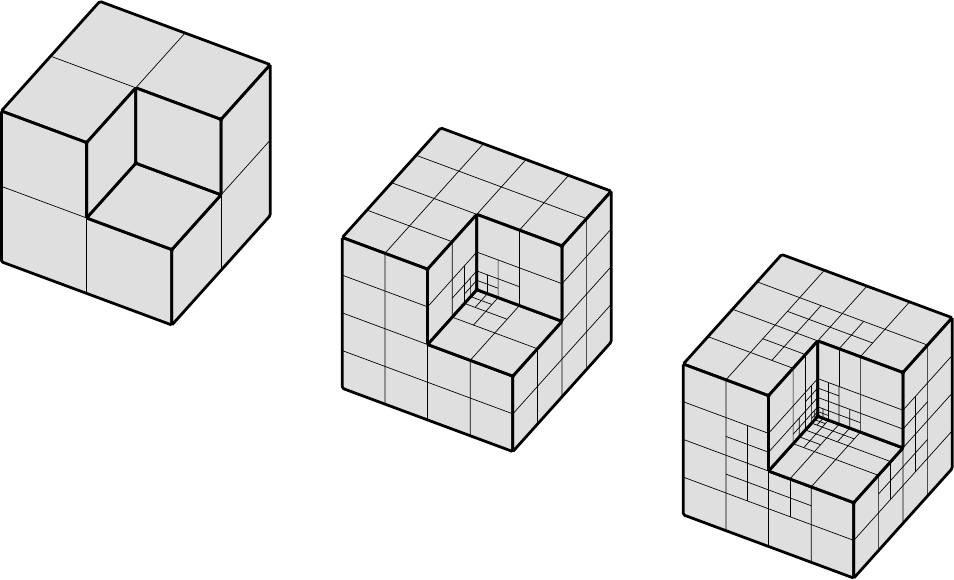}
The results of the isogeometric Nyström method are compared with an $h$-adaptive FEM solution, calculated by the Abaqus software package. 
The variation of the temperature $\primary$ is compared along a straight line defined by the endpoints $(0,0,1)$ and $(-0.5,-0.5,1)$ on $\Gamma_N$. The heat flux $\dual$ is compared along a straight line defined by $(0,0,0)$ and $(0.5,0.5,0)$ on $\Gamma_D$. Both lines are drawn in \figref{fig:fichera_BC} and the results are depicted in \figref{fig:fichera_results}.%
\tikzfig{tikz/ficheracube_results.tex}{1.0}{Temperature $\primary$ and heat flux $\dual$ along defined paths on the Fichera Cube calculated with different analysis methods}{fig:fichera_results}{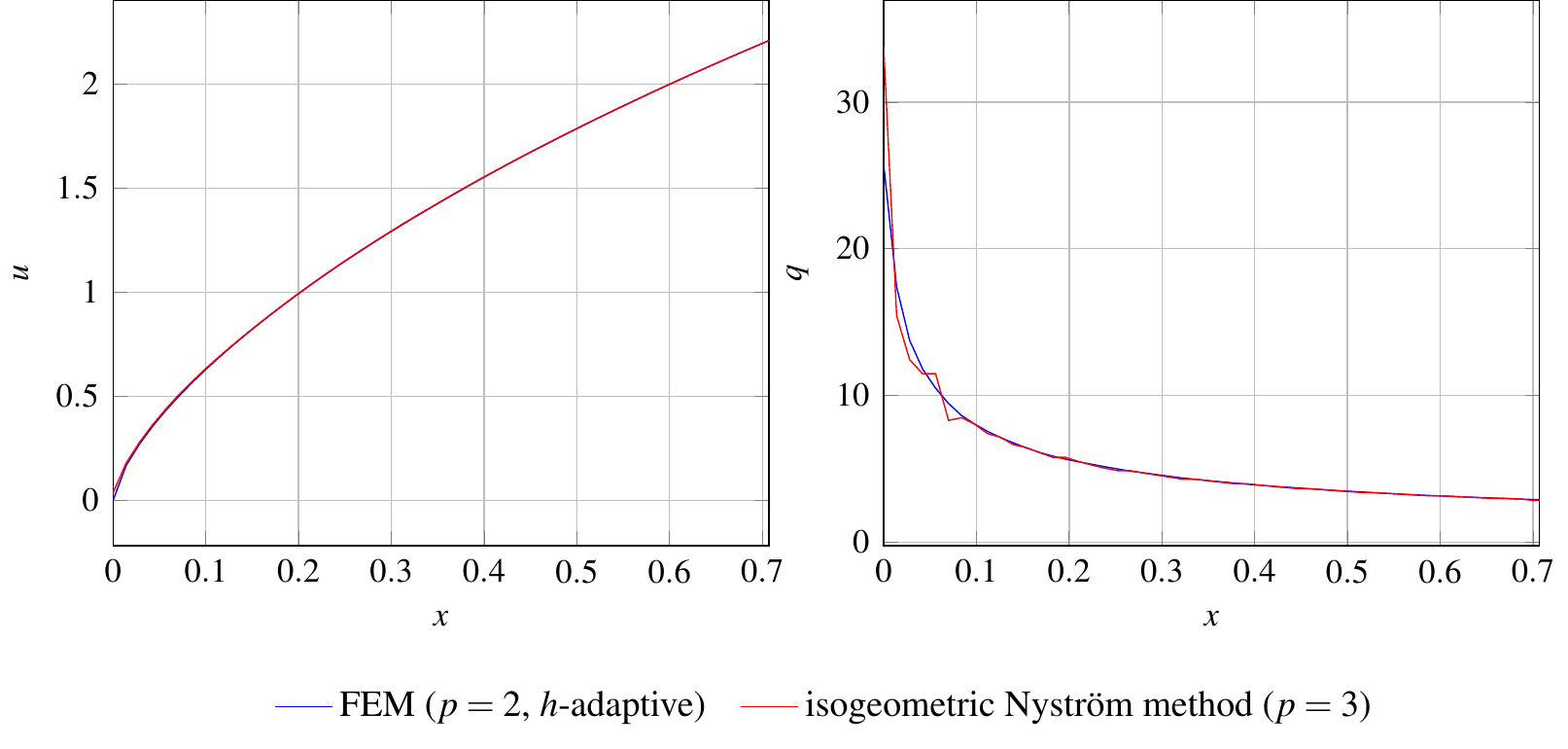}
The example on the Fichera cube demonstrates the ability of the isogeometric Nyström method to tackle physical singularities by using local refinement. 
\mycomment[GB]{Ausschlenker beim Heat-Flux erklären. Klarstellung, warum der Vergleich einer Singulären Lösung Sinn macht...}


%% file: results_spanner.tex
\subsubsection{Spanner}
\label{sec:res:spanner}

A single-ended open-jaw spanner is considered next. The geometry is defined by a CAGD model which is depicted on the top of \figref{fig:spanner}. The elastic constants are given by $E=\SI{2.1e5}{\mega\pascal}$ and $\nu=0.3$. Homogeneous Dirichlet boundary conditions, i.e. zero displacements $\pt{\primary}=0$ are applied at the open-jaw. The handle is subjected to a constant vertical loading $\vek{\dual}=( 0,0,t_z )$ with $t_z=\SI{15}{\mega\pascal}$. This setup and the partition of integration elements resulting from the CAGD model is illustrated in the middle of \figref{fig:spanner}.%
\tikzfig{tikz/spanner_all.tex}{0.15}{CAGD model of the spanner where thin black lines delimiting the NURBS patches (top). The partition into integration elements and the boundary conditions are depicted in the middle, the deflection as a result of linear elastic analysis is sketched on the bottom.}{fig:spanner}{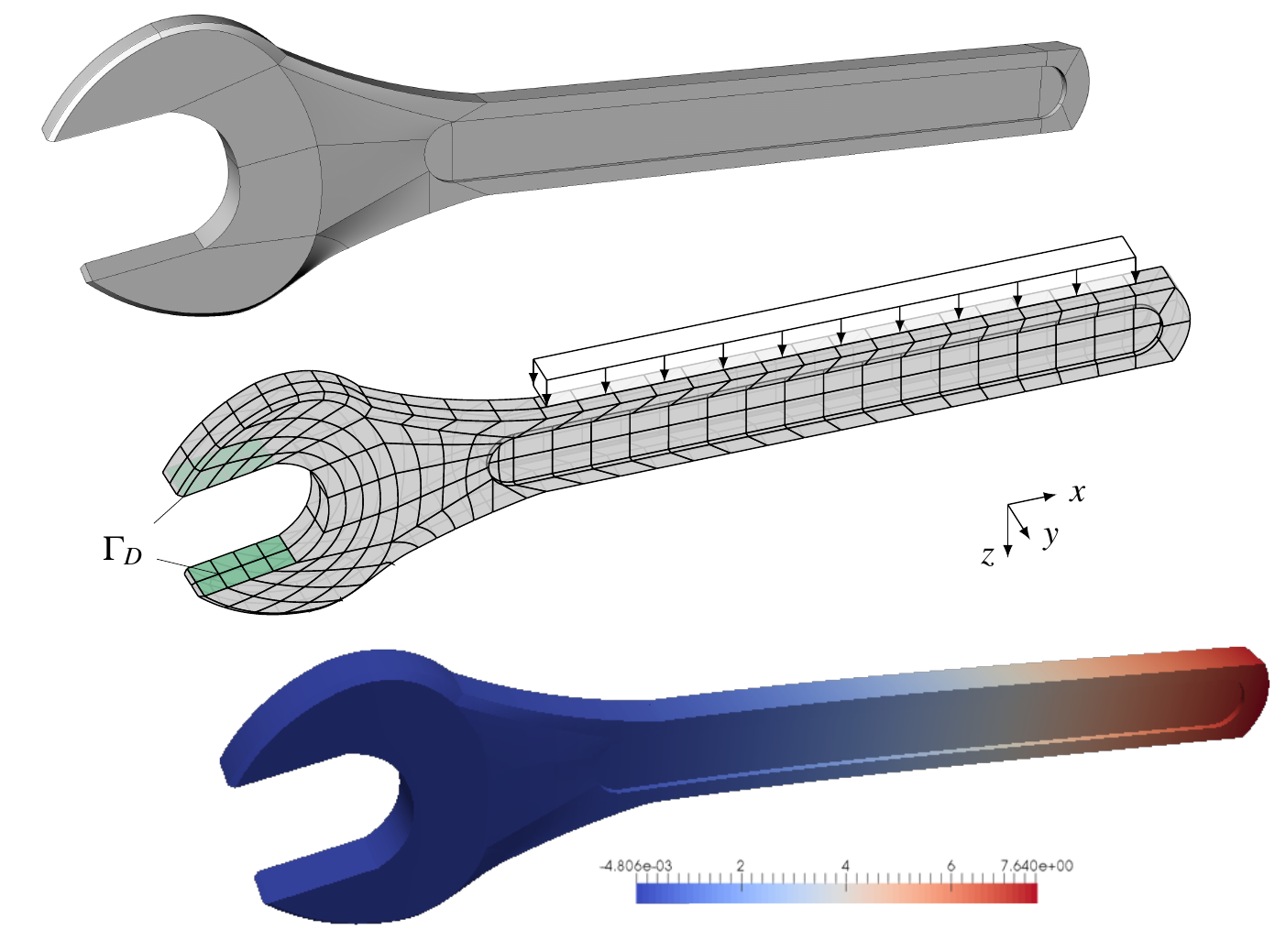}

The spanner is analyzed with the isogeometric Nyström method and the results post-processed with the strategy described in \secref{sec:iga:postprocessing}. The deflected tool is sketched at the bottom of \figref{fig:spanner}. The vertical and horizontal deflection at the top along a straight line is shown in \figref{fig:spanner_results}. The results are compared with an adaptive refined FEM analysis with Abaqus.%
\tikzfig{tikz/spanner_results.tex}{1.0}{Vertical deflection $u_z$ (left) and horizontal deflection $u_x$ (right) in [\si{\milli\metre}] at the top of the spanner along a straight line}{fig:spanner_results}{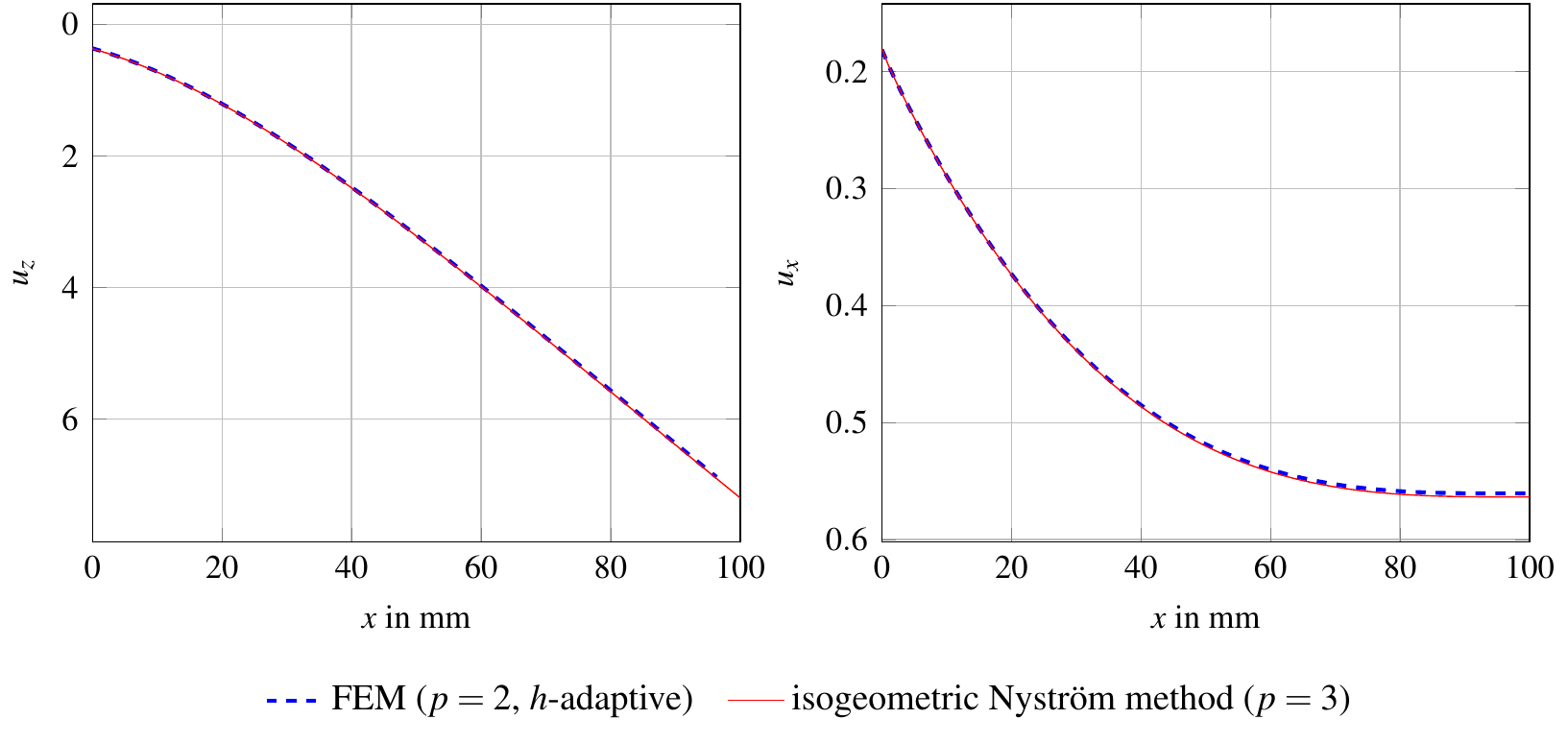}%

%% file: discussion.tex
\section{Discussion and Remarks}
\label{sec:discussion}

In this section, we further remark on the properties of the isogeometric Nyström method, practical considerations and on some implementation details.  
 
As for many other methods for solving BIEs, the accurate evaluation of singular and nearly singular integrals is essential. For locally corrected Nyström methods precise singular entries are achieved by means of properly evaluated entries for the right hand sides of the moment equation~\eqref{eq:LocalCorrectionSystem}. In the presented approach, the error of the integrals in equation~\eqref{eq:SL_LCN} is controllable and the integration is performed with an accuracy at least one magnitude lower than the lowest observed errors in the numerical tests of \secref{sec:results}. The quality of matrix entries related to nearly singular integration regions is controlled by the admissibility factor $\eta$. Depending on the complexity of the computational geometry, the factor is typically $\eta=\fromto{2}{6}$. However, increasing $\eta$ results to a larger locally corrected region and therefor to additional numerical effort for the analysis. A heuristic determination criterion for $\eta$ may assure the desired accuracy but an adaptive integration strategy like for BEM integral kernels is not viable due to the pointwise nature of the method. 

The presented results for heat conduction problems and elasticity in three dimensions demonstrate the practical applicability of the method. Although tensor product surfaces are used, local refinement is possible. The example on the Fichera cube demonstrates the ability of the isogeometric Nyström method to tackle physical singularities as well. The results for real world examples such as the cantilever beam or the spanner are of comparable quality to the standard FEM. This has been outlined in \secref{sec:res:fichera} and in \secref{sec:res:spanner}. For such practical problems, the $p$-refinement strategy is the first choice. However, for complicated geometries this approach is not always suitable. Hence, the convergence behavior for Fredholm integral equations \eqref{eq:FredholmFirst} of the first kind or direct BIEs like equation~\eqref{eq:DirectBIE} requires further attention to increase the robustness of the method. We would like to emphasize that the presented approach describes a straightforward implementation of the Nyström method with a strong focus on its implementation into a isogeometric framework. Implementations based on kernel splitting or other techniques may perform better, but such formulations depend on the underlying physical problem to be solved. In three dimensions efficient formulations are still a matter of research \cite{bremer2012a,bremer2014}.

Its pointwise nature greatly simplifies the implementation of the Nyström method into fast summation methods. The method has been considered in the landmark paper by \citet{rokhlin1985} for the fast multipole method. One of the advantages of the Nyström method compared to other methods solving boundary integral equations is, that point-wise supports are clearly related to spatial regions on which fast method are usually based upon. For instance, with Galerkin BEM supports of the basis functions may span multiple regions and therefor the partition and the treatment of near- and far-field system matrix entries gets more complex. The presented application makes use of hierarchical matrices \cite{hackbusch1999} 
adopting matrix operations of the HLib-library \cite{hlib}.

%% file: conclusion.tex
\section{Conclusion}
\label{sec:conclusion}

\mytodo[GB]{An isogeometric framework was applied to the Nystrom method. 
The Nystrom method presents an attractive alternative to classical BEM methods, its main feature being that no approximation of the unknown with basis functions is required.  The integral equations are solved by numerical integration where the unknown parameters are located at the Gauss points.}


In this work, an isogeometric framework was applied to the locally corrected Nyström method. The presented approach is suitable for any type of CAGD surface representation which provides a valid geometric mapping from parameter to real space. Hence, the method is applicable to NURBS, subdivision surfaces and T-splines in a straightforward way and can be easily adapted to further developments and new technologies in CAGD.

In this paper we explain the implementation of NURBS surfaces for the geometric description which are the commonly used in computer aided design. For such tensor product surfaces, the method inherently permits local $h$-refinement.

The isogeometric Nyström method presents an attractive alternative to classical methods solving boundary integral equations, its main feature being that no approximation of the unknown with basis functions is required. These equations are directly solved by numerical quadrature where the unknown parameters are located at the quadrature points.

The isogeometric Nyström method is characterized by a pointwise evaluation of the fundamental solution on the surface. Regularization of the singular integrals is carried out by means local correction with B-spline basis functions. 
These basis functions are also used for post-processing purposes, where the results at quadrature points are interpolated on the surface.

Due to its discrete, pointwise pattern, the Nyström method is well suited for fast summation methods such as the fast multipole method or hierarchical matrices. 


On test examples it has been shown that the method preforms well but that its convergence properties are different to commonly used methods such as the boundary element method. Some problems that need further attention were pointed out. It is hoped that this paper gives impetus for further investigation in this promising alternative to classical collocation type approaches.




%% file: appendixTransformation.tex
\section{Local Element Mapping}
\label{appsec:appendixTransformation}

The mapping from the initial element $\iegLocal_{0}$ defined by a tensor product of $\KVC_\totalA$ and $\KVC_\totalB$ to local elements $\iegLocal_{\ell}$ specified by refinement points $\locRefInfo_\indexLevel$ is discussed.
In a first step, each $\iegLocal_{0}$ is represented by means of its corner nodes, i.e. the knots of its corresponding knots span $\vek{\indexSpan}=(\indexSpan_1,\indexSpan_{2})^\trans$.
They are summarized in a \emph{node matrix} $\nodeMatrix_0$ such that
\begin{align}	
	\nodeMatrix_0 = 	
	\begin{pmatrix}
		\KVC_\totalA[\indexSpan_1] & \KVC_\totalA[\indexSpan_1+1] & \KVC_\totalA[\indexSpan_1] & \KVC_\totalA[\indexSpan_1+1] \\
		\KVC_\totalB[\indexSpan_2] & \KVC_\totalB[\indexSpan_2+1] & \KVC_\totalB[\indexSpan_2+1] & \KVC_\totalB[\indexSpan_2] \\
		1 & 1 & 1 & 1
	\end{pmatrix}.
\end{align}
The mapping from $\iegLocal_{0}$ to its $\iegLocal_{1,\indexA}$ includes the translation and scaling of the corner nodes.
They are assembled in a \emph{transformation matrix} $\transMatrix_{1,\indexA}$ which is defined for each $\iegLocal_{1,\indexA}$ by
\begin{align}	
	\transMatrix_{1,\indexA} = 	
	\begin{pmatrix}
		l_{1,\uu_1}/ l_{0,{\uu_1}} & 0 & t_{\uu_1} \\
		0 &  l_{1,\uu_2} / l_{0,\uu_2}  & t_{\uu_2}  \\
		0 & 0 & 1 
	\end{pmatrix}.
\end{align}
The last column refers to the translation $t$ of the first corner node, whereas the diagonal entries are related to the lengths of the initial ($l_0$) and refined element ($l_1$)  in each parametric direction~$\uu$.
The construction of $\transMatrix_{1,\indexA}$ due to a set of refinement points $\locRefInfo$ of the first level is summarized in \algref{alg:Transformation}.
\begin{myalgorithm}{Set up transformation matrices for the next level}{alg:Transformation}
	\REQUIRE Node matrix $\nodeMatrix$ and related refinement points $\locRefInfo$ of the subsequent level $\ell$
	\STATE $\Lieg_{\uu_1} =  \nodeMatrix[1,3] - \nodeMatrix[1,1] $ 
	\STATE $\Lieg_{\uu_2} =  \nodeMatrix[2,3] - \nodeMatrix[2,1] $ 
	\STATE initialize temporary knot vectors $\KVC_\totalA$ and $\KVC_\totalB$
	\STATE $\KVC_\totalA \myinsert \nodeMatrix[1,k], k=1,3$
	\STATE $\KVC_\totalB \myinsert \nodeMatrix[2,k], k=1,3$
	\FORALL{$\locRefInfo$}
		\STATE $\KVC_\totalA \myinsert \locRefInfo_\indexA[1]$
		\STATE $\KVC_\totalB \myinsert \locRefInfo_\indexA[2]$
	\ENDFOR
	\STATE $\LiegLocal_\totalA \store$ length of each non-zero knot spans of $\KVC_\totalA$
	\STATE $\LiegLocal_\totalB \store$ length of each non-zero knot spans of $\KVC_\totalB$
	\STATE initialize array $\vek{a}_T$ for transformation matrices $\transMatrix$
	\STATE $t_{\uu_2} \store 0$  
	\FORALL{$\LiegLocal_\indexB \in \LiegLocal_\totalB$}
		\STATE $t_{\uu_1} \store 0$  
		\FORALL{$\LiegLocal_\indexA \in \LiegLocal_\totalA$}
			\STATE $\transMatrix = \diag(\LiegLocal_\indexA /  \Lieg_{\uu_1}, \LiegLocal_\indexB / \Lieg_{\uu_2},1)$
			\STATE $\transMatrix[1,3] = \nodeMatrix[1,1] \left(1-\LiegLocal_\indexA /  \Lieg_{\uu_1} \right) + t_{\uu_1}$
			\STATE $\transMatrix[2,3] = \nodeMatrix[2,1] \left(1-\LiegLocal_\indexB / \Lieg_{\uu_2} \right) + t_{\uu_2}$
			\STATE $\vek{a}_T \myinsert \transMatrix$
			\STATE $t_{\uu_1} \store t_{\uu_1} + \LiegLocal_\indexA$  
		\ENDFOR
		\STATE $t_{\uu_2} \store t_{\uu_2} + \LiegLocal_\indexB$  
	\ENDFOR
	\RETURN $\vek{a}_T$
\end{myalgorithm}
Since the nodes of $\nodeMatrix_0$ are represented in homogeneous coordinates with $\w = 1$, the transformation to the nodes $\nodeMatrix_{1,\indexA}$ of local elements $\iegLocal_{1,\indexA}$ can be expressed by a matrix product as 
\begin{align}
	\label{eq:defLocalElement}
	\iegLocal_{1,\indexA} \coloneqq \nodeMatrix_{1,\indexA} = \transMatrix_{1,\indexA} \:\nodeMatrix_0.
\end{align}
If there are refinement points $\locRefInfo_\ell$ of a higher level, i.e. $\ell > 1$, within an $\iegLocal_{1,\indexA}$ additional transformation matrices $\transMatrix_{2,\indexA}$ are constructed based on $\nodeMatrix_{1,\indexA}$ and $\locRefInfo_\ell$.
The resulting local elements $\iegLocal_{2,\indexA}$ are given by
\begin{align}
	\label{eq:defLocalElementLevel2}
	\iegLocal_{2,\indexA} \coloneqq \nodeMatrix_{2,\indexA} = \transMatrix_{2,\indexA} \:\nodeMatrix_{1,\indexB} = \transMatrix_{2,\indexA} \: \transMatrix_{1,\indexB} \:\nodeMatrix_0.
\end{align}
The accumulated transformation matrices $\hat{\transMatrix}_{\indexLevel,\indexA}$ relate the final $\iegLocal_{\indexLevel,\indexA}$ due to all refinement level to the initial knot span $\nodeMatrix_0$.
\begin{align}
	\label{eq:defLocalElementMoreLevels}
	\iegLocal_{\indexLevel,\indexA} & \coloneqq \nodeMatrix_{\indexLevel,\indexA} = \hat{\transMatrix}_{\indexLevel,\indexA} \: \nodeMatrix_0 
	&& \with && \hat{\transMatrix}_{\indexLevel,\indexA} = \prod_{\indexC\in\totalLevel} \transMatrix_{\indexC,\indexD}
\end{align}
where $\totalLevel$ denotes an index set of all levels defining $\iegLocal_{\indexLevel,\indexA}$ which is ordered decreasingly.
The set up of $\hat{\transMatrix}_{\indexLevel,\indexA}$ is described in \algref{alg:moreLevels}. 
\begin{myalgorithm}{Hierarchical refinement}{alg:moreLevels}
	\REQUIRE Node matrix $\nodeMatrix_0$ of an element $\ieg$ and refinement points $\locRefInfo$ of all levels
	\STATE initialize array $\vek{a}_T$ for transformation matrices $\hat{\transMatrix}_\indexLevel \in \R^{3 \times 3}$
	\STATE $\vek{a}_T \myinsert \transMatrix_0 = \diag\left(1,1,1\right)$ 
	\FORALL{refinement levels $\indexLevel$}
		\STATE initialize temporary array $\vek{b}_T$ for $\transMatrix_{\indexLevel}$
		\FORALL[loop over all computed local elements $\iegLocal_{\indexLevel-1}$]{$\transMatrix_{\indexLevel-1,\indexC} \in \vek{a}_T$}
			\STATE $\nodeMatrix_{\indexLevel-1,\indexC}$ \store~$\transMatrix_{\indexLevel-1,\indexC} \nodeMatrix_0$
			\STATE initialize temporary array $\vek{c}_{\locRefInfo}$ for refinement points
	;		\FORALL{$\locRefInfo_{\indexLevel,\indexD} \in \locRefInfo_{\indexLevel}$}
				\IF{$\locRefInfo_{\indexLevel,\indexD}$ is inside $\iegLocal_{\indexLevel-1,\indexC}$ related to $\nodeMatrix_{\indexLevel-1,\indexC}$}
					\STATE $\vek{c}_{\locRefInfo} \myinsert \locRefInfo_{\indexLevel,\indexD}$
				\ENDIF
			\ENDFOR
			\IF{$\vek{c}_{\locRefInfo} = \emptyset $}
				\STATE $\vek{b}_t \myinsert \transMatrix_{\indexLevel-1,\indexC}$ 
			\ELSE
				\STATE $\vek{c}_T$ \store~array of $\transMatrix_\indexLevel$ set up by \algref{alg:Transformation} with $\nodeMatrix_{\indexLevel-1,\indexC}$ and $\vek{c}_{\locRefInfo}$
				\FORALL{$\transMatrix_{\indexLevel,r} \in \vek{c}_T$}
					\STATE $\vek{b}_T \myinsert \transMatrix_{\indexLevel,r} \: \transMatrix_{\indexLevel-1,\indexC}  $ 
				\ENDFOR
			\ENDIF
		\ENDFOR
		\STATE $\vek{a}_T \store \vek{b}_T$
	\ENDFOR
	\RETURN $\vek{a}_T$
\end{myalgorithm}


%% file: literature.tex
\clearpage

\bibliographystyle{elsarticle-num-names}

\bibliography{mybibliography}


